\documentclass[12pt]{amsart}
\usepackage[utf8]{inputenc}
\usepackage[french]{babel}
\usepackage{geometry} 
\usepackage{graphicx}				
\usepackage{amssymb}
\usepackage{graphicx}
\usepackage{epstopdf}
\usepackage{color}
\usepackage{pgf,tikz}
\usetikzlibrary{arrows}
\usepackage[colorlinks=true,citecolor=cyan,urlcolor=blue]{hyperref} 
   
\newtheoremstyle{mes_theoremes}{1.5em}{1.9em}{}{}{\bfseries}{~:~}{\parskip}{\thmname{#1}\thmnumber{ #2}\thmnote{ (#3)}}

\theoremstyle{mes_theoremes}

\newtheorem*{pre}{Preuve}
\newtheorem{note}{Notation}
\newtheorem{de}{Définition}

\newtheorem{prop}{Proposition}

\newtheorem*{rmq}{Remarque}

\newtheorem*{ex*}{Exemple}

\newtheoremstyle{mes_preuves}{1.5em}{2em}{}{}{\it}{~:~}{\parskip}{\thmname{#1}\thmnumber{}}

\theoremstyle{mes_preuves}

\newcommand{\sch}{\mathop{\mathrm{Sch}}}
\newcommand{\s}{\mathop{\mathrm{S}}}
\newcommand{\dd}{\mathop{\mathrm{Dyck}}}
\newcommand{\ag}{\mathop{\mathrm{aireg}}}
\newcommand{\A}{\mathop{\mathrm{aire}}}
\newcommand{\dinv}{\mathop{\mathrm{dinv}}}
\newcommand{\p}{\mathop{\mathbb{P}}}
\newcommand{\PP}{\mathop{\mathrm{G}}}
\newcommand{\PPE}{\mathop{\mathrm{E}}}
\newcommand{\Pag}{\mathop{\mathrm{F}}}
\newcommand{\hexa}{\mathop{\mathrm{H}}}
\newcommand{\Ps}{\mathop{\mathrm{Psch}}}
\newcommand{\Cs}{\mathop{\mathrm{C_hsch}}}{}
\newcommand{\F}{\infty \hspace{-9pt}\begin{Huge}^\text{O}_\text{O}\end{Huge} \hspace{-2pt}_\smile \hspace{-14pt}\bullet}

\title{Les chemins de Schröder
\\(Rapport de stage, été 2014)}

\author{\href{mailto:wallace.nancy@courrier.uqam.ca}{Nancy Wallace}} 
\begin{document}
\maketitle

J'aimerais d'abord remercier François Bergeron pour ses communications personnelles. Elles sont la référence principale de ce document. Les références seront notées seulement si ce n'est pas le cas.
\\

Après avoir posé les définitions nécessaires à la compréhension du sujet, nous discuterons de statistique d'inversion diagonale dans les $r$-Schröder, de chemins de stationnement dans les $r$-Schröder à pente entière et nous développerons une formule pour les chemins de Schröder ayant une fraction unitaire comme pente.
\tableofcontents
\newpage
\section{Définitions}
\hspace{20pt}
\begin{de} Un \begin{bf} pas\end{bf}, est un déplacement d'une longueur prédéterminée, d'un point à coordonnées entières vers à un autre point à coordonnées entières, dans le plan cartésien. \end{de}

\begin{de} Un \begin{bf} chemin de Schröder\end{bf} est une suite de pas débutant en $(0,n)$ et se terminant en $(n,0)$ sans jamais passer au-dessus de la droite passant par $(0,n)$ et $(n,0)$. De plus, les seuls pas autorisé sont $(0,-1)$, dit vers le bas (certains auteurs disent aussi vers le sud), $(1,0)$ dits vers la droite (vers l'est) ou $(1,-1)$ un pas diagonal. Ceci nous amène naturellement à la définition suivante.\end{de}
 
 \begin{de} Un \begin{bf} chemin r-Schröder \end{bf} est une suite de pas débutant en $(0,n)$ et se terminant en $(rn,0)$ sans jamais passer au-dessus de la droite passant par $(0,n)$ et $(rn,0)$. De plus, les seuls pas autorisé sont $(0,-1)$, dit vers le bas, $(1,0)$ dits vers la droite, ou $(r,-1)$ un pas diagonal. \end{de}
 
 Ces définitions impliquent clairement que le chemin est entièrement compris dans le triangle $(0,0), (0,n), (rn,0)$. De plus, avec cette définition on remarque qu'un chemin de Dyck est un chemin de Schröder n'ayant pas de pas diagonal.

 \begin{note} Notons $\mathbf{\sch}_{n,d}$, l'ensemble des chemins de Schröder ayant $n-d$ pas diagonaux et $\mathbf{\sch}_{n,d}^r$ , l'ensemble des chemins $r$-Schröder avec $n-d$ pas diagonaux. L'ensemble des $r$-Dyck est noté $\mathbf{\dd}_n^{r}$, donc $\mathbf{\dd}_n^{r}=\mathbf{\sch}_{n,n}^r$. Enfin, l'ensemble des chemins $r$-Schröder est noté $\mathbf{\sch}_{n}^r$. \end{note}

  De plus, la formule suivante donne le nombre de chemins, ayant $n-d$ pas diagonaux: 
\[ |\sch {_{n,d}^r}|=\s{_{n,d}^{ r}}=\frac{1}{dr+1}\binom{n}{d}\binom{dr+n}{n} .\]

Remarquons que pour $d=n$ ceci correspond au nombre de Fuss-Catalan donné par:
\[ |\dd{_n^{~r}}| =\frac{1}{nr+1}\binom{nr+n}{n} =|\sch{_{n,n}^r}|. \]
  \begin{de} Une ligne $i$ est comprise entre la droite $i$ et la droite $i+1$.  Ici les droites sont numérotées du haut vers le bas de façon à ce que la première soit la droite horizontale passant par $(0,n)$ et la droite $n+1$ est l'axe des $X$. Soit $\alpha \in \sch_{n,d}^r$, \begin{bf} l'aire gauche d'une ligne\end{bf}, notée $\mathbf{\ag}_i(\alpha)$, est le nombre de carrés pleins se trouvant entre l'axe des $Y$ et le chemin $\alpha$. \end{de}
 
\begin{note}  Soit $\alpha \in \sch_{n,d}^r$ . Pour une rangée $i$ ayant une aire gauche de valeur $k$, notons $k$ l'aire gauche de cette rangée si le chemin entre la ligne $i-1$ et la ligne $i$ est un pas vers le bas et notons $\bar{k}$ l'aire gauche de cette rangée si le chemin $\alpha$ entre la ligne $i$ et la ligne $i+1$ est un pas diagonal. Par abus de notation, la suite de ces aires est notée $\alpha$, car cela donne un codage bijectif des chemins de Schröder et nous pouvons identifier de façon unique $\alpha$ à ce codage.\\  \end{note}
\vspace{-30pt}  
\definecolor{ffcqcb}{rgb}{1.0,0.7529411764705882,0.796078431372549}
\definecolor{yqqqqq}{rgb}{0.5019607843137255,0.0,0.0}
\definecolor{xfqqff}{rgb}{0.4980392156862745,0.0,1.0}
\definecolor{qqqqff}{rgb}{0.0,0.0,1.0}
\definecolor{xdxdff}{rgb}{0.49019607843137253,0.49019607843137253,1.0}
\definecolor{uuuuuu}{rgb}{0.26666666666666666,0.26666666666666666,0.26666666666666666}
\begin{tikzpicture}[line cap=round,line join=round,>=triangle 45,x=.7cm,y=.7cm]
\clip(-3.759999999999999,-2.2000000000000033) rectangle (10.440000000000007,6.08);
\fill[color=xfqqff,fill=xfqqff,fill opacity=0.1] (0.0,-0.0) -- (0.0,5.0) -- (10.0,0.0) -- cycle;
\fill[color=ffcqcb,fill=ffcqcb,fill opacity=0.45] (0.0,3.0) -- (2.0,3.0) -- (2.0,2.0) -- (0.0,2.0) -- cycle;
\fill[color=ffcqcb,fill=ffcqcb,fill opacity=0.45] (0.0,2.0) -- (0.0,1.0) -- (3.0,1.0) -- (3.0,2.0) -- cycle;
\fill[color=ffcqcb,fill=ffcqcb,fill opacity=0.45] (0.0,1.0) -- (0.0,-0.0) -- (5.0,0.0) -- (5.0,1.0) -- cycle;
\draw [color=xfqqff] (0.0,-0.0)-- (0.0,5.0);
\draw [color=xfqqff] (0.0,5.0)-- (10.0,0.0);
\draw [color=xfqqff] (10.0,0.0)-- (0.0,-0.0);
\draw [line width=5.2pt,color=yqqqqq] (0.0,5.0)-- (0.0,4.0);
\draw [line width=5.2pt,color=yqqqqq] (0.0,4.0)-- (2.0,3.0);
\draw [line width=5.2pt,color=yqqqqq] (2.0,3.0)-- (2.0,2.0);
\draw [line width=5.2pt,color=yqqqqq] (2.0,2.0)-- (3.0,2.0);
\draw [line width=5.2pt,color=yqqqqq] (3.0,2.0)-- (3.0,1.0);
\draw [line width=5.2pt,color=yqqqqq] (3.0,1.0)-- (5.0,1.0);
\draw [line width=5.2pt,color=yqqqqq] (5.0,1.0)-- (7.0,0.0);
\draw [line width=5.2pt,color=yqqqqq] (7.0,0.0)-- (10.0,0.0);
\draw [color=ffcqcb] (0.0,3.0)-- (2.0,3.0);
\draw [color=ffcqcb] (2.0,3.0)-- (2.0,2.0);
\draw [color=ffcqcb] (2.0,2.0)-- (0.0,2.0);
\draw [color=ffcqcb] (0.0,2.0)-- (0.0,3.0);
\draw [color=ffcqcb] (0.0,2.0)-- (0.0,1.0);
\draw [color=ffcqcb] (0.0,1.0)-- (3.0,1.0);
\draw [color=ffcqcb] (3.0,1.0)-- (3.0,2.0);
\draw [color=ffcqcb] (3.0,2.0)-- (0.0,2.0);
\draw [color=ffcqcb] (0.0,1.0)-- (0.0,-0.0);
\draw [color=ffcqcb] (0.0,-0.0)-- (5.0,0.0);
\draw [color=ffcqcb] (5.0,0.0)-- (5.0,1.0);
\draw [color=ffcqcb] (5.0,1.0)-- (0.0,1.0);
\draw (-3.6799999999999984,5.800000000000001) node[anchor=north west] {\footnotesize{Ligne}};
\draw (-3.3599999999999985,4.86) node[anchor=north west] {$1$};
\draw (-3.3599999999999985,3.8800000000000003) node[anchor=north west] {$2$};
\draw (-3.3599999999999985,2.9) node[anchor=north west] {$3$};
\draw (-3.3599999999999985,1.8999999999999995) node[anchor=north west] {$4$};
\draw (-3.3599999999999985,0.8999999999999992) node[anchor=north west] {$5$};
\draw [color=xfqqff] (0.0,4.0)-- (2.0,4.0);
\draw [color=xfqqff] (0.0,3.0)-- (4.0,3.0);
\draw [color=xfqqff] (0.0,2.0)-- (6.0,2.0);
\draw [color=xfqqff] (0.0,1.0)-- (8.0,1.0);
\draw [color=xfqqff] (1.0319999999999998,4.484)-- (1.0,0.0);
\draw [color=xfqqff] (2.0,4.0)-- (2.0,0.0);
\draw [color=xfqqff] (3.0239999999999996,3.4880000000000004)-- (3.0,0.0);
\draw [color=xfqqff] (4.0,3.0)-- (4.0,0.0);
\draw [color=xfqqff] (5.032,2.484)-- (5.0,0.0);
\draw [color=xfqqff] (6.0,2.0)-- (6.0,0.0);
\draw [color=xfqqff] (7.008,1.496)-- (7.0,0.0);
\draw [color=xfqqff] (8.0,1.0)-- (8.0,0.0);
\draw [color=xfqqff] (9.008,0.49600000000000044)-- (9.0,0.0);
\draw (-2.119999999999998,5.800000000000001) node[anchor=north west] {\footnotesize{Aire gauche}};
\draw (-1.3399999999999976,4.9) node[anchor=north west] {$0$};
\draw (-1.3399999999999976,3.9000000000000004) node[anchor=north west] {$\bar{0}$};
\draw (-1.3399999999999976,2.9) node[anchor=north west] {$2$};
\draw (-1.3399999999999976,1.9199999999999995) node[anchor=north west] {$3$};
\draw (-1.3399999999999976,0.9199999999999992) node[anchor=north west] {$\bar{5}$};
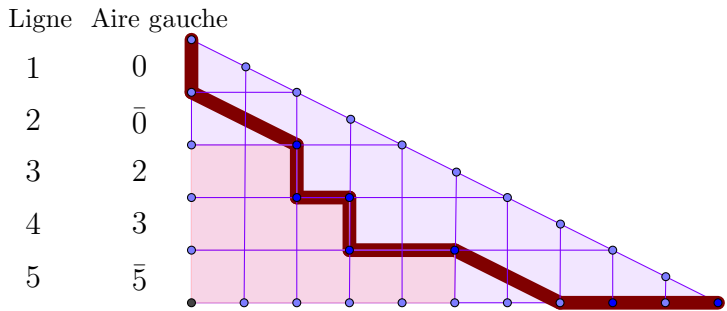
\draw (-0.5199999999999982,-0.8000000000000014) node[anchor=north west] {\footnotesize{$\text{Chemin de Schröder } \alpha \in Sch_{5,2}^2 \text{ de codage } 0\bar{0}23\bar{5}$}};
\begin{scriptsize}
\draw [fill=uuuuuu] (0.0,-0.0) circle (1.5pt);
\draw [fill=xdxdff] (0.0,5.0) circle (1.5pt);
\draw [fill=qqqqff] (10.0,0.0) circle (1.5pt);
\draw [fill=xdxdff] (0.0,4.0) circle (1.5pt);
\draw [fill=qqqqff] (2.0,3.0) circle (1.5pt);
\draw [fill=qqqqff] (2.0,2.0) circle (1.5pt);
\draw [fill=qqqqff] (3.0,2.0) circle (1.5pt);
\draw [fill=qqqqff] (3.0,1.0) circle (1.5pt);
\draw [fill=qqqqff] (5.0,1.0) circle (1.5pt);
\draw [fill=xdxdff] (7.0,0.0) circle (1.5pt);
\draw [fill=xdxdff] (0.0,3.0) circle (1.5pt);
\draw [fill=xdxdff] (0.0,2.0) circle (1.5pt);
\draw [fill=xdxdff] (0.0,1.0) circle (1.5pt);
\draw [fill=xdxdff] (5.0,0.0) circle (1.5pt);
\draw [fill=xdxdff] (2.0,4.0) circle (1.5pt);
\draw [fill=xdxdff] (4.0,3.0) circle (1.5pt);
\draw [fill=xdxdff] (6.0,2.0) circle (1.5pt);
\draw [fill=xdxdff] (8.0,1.0) circle (1.5pt);
\draw [fill=xdxdff] (1.0319999999999998,4.484) circle (1.5pt);
\draw [fill=xdxdff] (1.0,0.0) circle (1.5pt);
\draw [fill=xdxdff] (2.0,0.0) circle (1.5pt);
\draw [fill=xdxdff] (3.0239999999999996,3.4880000000000004) circle (1.5pt);
\draw [fill=xdxdff] (3.0,0.0) circle (1.5pt);
\draw [fill=xdxdff] (4.0,0.0) circle (1.5pt);
\draw [fill=xdxdff] (5.032,2.484) circle (1.5pt);
\draw [fill=xdxdff] (6.0,0.0) circle (1.5pt);
\draw [fill=xdxdff] (7.008,1.496) circle (1.5pt);
\draw [fill=qqqqff] (8.0,0.0) circle (1.5pt);
\draw [fill=xdxdff] (9.008,0.49600000000000044) circle (1.5pt);
\draw [fill=xdxdff] (9.0,0.0) circle (1.5pt);
\end{scriptsize}
\end{tikzpicture}

 \begin{de} Pour $\alpha \in \sch_{n,d}^r$, l'\begin{bf} aire de la ligne $i$, \end{bf} notée $\mathbf{\A}_i(\alpha)$ est le nombre de pas vers la droite situé sur la $i$-ème ligne entre le chemin $\alpha$ et la droite passants par $(0,n)$ et $(rn,0)$.  \end{de}
 
  \begin{de} Pour $\alpha \in \sch_{n,d}^r$ ,  l'\begin{bf}aire de $\alpha$\end{bf} noté $\mathbf{\A}(\alpha)$ est la somme:
  \[ \A(\alpha)=\sum_{i=1}^{n} \A{_i}(\alpha).\] 
On a pour cet air la formule:
\[  \A(\alpha)=r(i-1)-|a_i|, \text{ où }, a_i=\ag{_i}(\alpha), |a_i|=a_i \text{ et } |\bar{a_i}|=a_i.\]  \end{de}

\begin{tikzpicture}[line cap=round,line join=round,>=triangle 45,x=1.0cm,y=1.0cm]
\clip(-1.6799999999999984,-2.380000000000001) rectangle (18.18000000000001,5.900000000000002);
\fill[color=xfqqff,fill=xfqqff,fill opacity=0.1] (0.0,-0.0) -- (0.0,5.0) -- (10.0,0.0) -- cycle;
\draw [color=xfqqff] (0.0,-0.0)-- (0.0,5.0);
\draw [color=xfqqff] (0.0,5.0)-- (10.0,0.0);
\draw [color=xfqqff] (10.0,0.0)-- (0.0,-0.0);
\draw (-1.319999999999998,5.640000000000002) node[anchor=north west] {Ligne};
\draw (-1.099999999999998,4.860000000000002) node[anchor=north west] {$1$};
\draw (-1.099999999999998,3.8800000000000012) node[anchor=north west] {$2$};
\draw (-1.099999999999998,2.9000000000000012) node[anchor=north west] {$3$};
\draw (-1.099999999999998,1.9000000000000008) node[anchor=north west] {$4$};
\draw (-1.099999999999998,0.9000000000000004) node[anchor=north west] {$5$};
\draw [color=xfqqff] (0.0,4.0)-- (2.0,4.0);
\draw [color=xfqqff] (0.0,3.0)-- (4.016,2.992);
\draw [color=xfqqff] (0.0,2.0)-- (6.0,2.0);
\draw [color=xfqqff] (0.0,1.0)-- (8.0,1.0);
\draw [color=xfqqff] (1.0319999999999998,4.484)-- (1.0,0.0);
\draw [color=xfqqff] (2.0,4.0)-- (2.0,0.0);
\draw [color=xfqqff] (3.0239999999999996,3.4880000000000004)-- (3.0,0.0);
\draw [color=xfqqff] (4.016,2.992)-- (4.0,0.0);
\draw [color=xfqqff] (5.032,2.484)-- (5.0,0.0);
\draw [color=xfqqff] (6.0,2.0)-- (6.0,0.0);
\draw [color=xfqqff] (7.008,1.496)-- (7.0,0.0);
\draw [color=xfqqff] (8.0,1.0)-- (8.0,0.0);
\draw [color=xfqqff] (9.008,0.49600000000000044)-- (9.0,0.0);
\draw (10.480000000000008,5.580000000000002) node[anchor=north west] {Aire };
\draw (10.620000000000008,4.880000000000002) node[anchor=north west] {$0$};
\draw (10.620000000000008,3.8800000000000012) node[anchor=north west] {$1$};
\draw (10.620000000000008,2.8800000000000012) node[anchor=north west] {$3$};
\draw (10.620000000000008,1.9000000000000008) node[anchor=north west] {$0$};
\draw (10.620000000000008,0.9000000000000004) node[anchor=north west] {$0$};
\draw (0.06000000000000259,-0.8200000000000003) node[anchor=north west] {$\text{Chemin de Schröder } \alpha \in Sch_{5,1}^2 \ \text{ d'aire  } aire(\alpha)=0+1+3+0+0=4$};
\draw [line width=5.2pt,color=yqqqqq] (0.0,5.0)-- (0.0,4.0);
\draw [line width=5.2pt,color=yqqqqq] (0.0,4.0)-- (1.0285459411239963,4.0);
\draw [line width=5.2pt,color=yqqqqq] (1.0285459411239963,4.0)-- (1.0213949355841847,2.997965348733896);
\draw [line width=5.2pt,color=yqqqqq] (6.0,2.0)-- (8.0,1.0);
\draw [line width=5.2pt,color=yqqqqq] (8.0,1.0)-- (8.0,0.0);
\draw [line width=5.2pt,color=yqqqqq] (8.0,0.0)-- (10.0,0.0);
\draw [line width=5.2pt,color=ffcqcb] (1.0285459411239963,4.0)-- (2.0,4.0);
\draw [line width=5.2pt,color=ffcqcb] (1.0213949355841847,2.997965348733896)-- (4.016,2.992);
\draw [line width=5.2pt,color=yqqqqq] (6.0,2.0)-- (3.0137614678899083,2.0);
\draw [line width=5.2pt,color=yqqqqq] (1.0213949355841847,2.997965348733896)-- (3.0137614678899083,2.0);
\begin{scriptsize}
\draw [fill=uuuuuu] (0.0,-0.0) circle (1.5pt);
\draw [fill=xdxdff] (0.0,5.0) circle (1.5pt);
\draw [fill=qqqqff] (10.0,0.0) circle (1.5pt);
\draw [fill=xdxdff] (0.0,4.0) circle (1.5pt);
\draw [fill=qqqqff] (4.0,3.0) circle (1.5pt);
\draw [fill=qqqqff] (5.98,1.98) circle (1.5pt);
\draw [fill=xdxdff] (7.0,0.0) circle (1.5pt);
\draw [fill=xdxdff] (0.0,3.0) circle (1.5pt);
\draw [fill=xdxdff] (0.0,2.0) circle (1.5pt);
\draw [fill=xdxdff] (0.0,1.0) circle (1.5pt);
\draw [fill=xdxdff] (5.0,0.0) circle (1.5pt);
\draw [fill=xdxdff] (2.0,4.0) circle (1.5pt);
\draw [fill=xdxdff] (4.016,2.992) circle (1.5pt);
\draw [fill=xdxdff] (6.0,2.0) circle (1.5pt);
\draw [fill=xdxdff] (8.0,1.0) circle (1.5pt);
\draw [fill=xdxdff] (1.0319999999999998,4.484) circle (1.5pt);
\draw [fill=xdxdff] (1.0,0.0) circle (1.5pt);
\draw [fill=xdxdff] (2.0,0.0) circle (1.5pt);
\draw [fill=xdxdff] (3.0239999999999996,3.4880000000000004) circle (1.5pt);
\draw [fill=xdxdff] (3.0,0.0) circle (1.5pt);
\draw [fill=xdxdff] (4.0,0.0) circle (1.5pt);
\draw [fill=xdxdff] (5.032,2.484) circle (1.5pt);
\draw [fill=xdxdff] (6.0,0.0) circle (1.5pt);
\draw [fill=xdxdff] (7.008,1.496) circle (1.5pt);
\draw [fill=qqqqff] (8.0,0.0) circle (1.5pt);
\draw [fill=xdxdff] (9.008,0.49600000000000044) circle (1.5pt);
\draw [fill=xdxdff] (9.0,0.0) circle (1.5pt);
\draw [fill=uuuuuu] (1.0285459411239963,4.0) circle (1.5pt);
\draw [fill=uuuuuu] (1.0213949355841847,2.997965348733896) circle (1.5pt);
\draw [fill=uuuuuu] (3.0137614678899083,2.0) circle (1.5pt);
\end{scriptsize}
\end{tikzpicture}
 
\section{Statistique d'inversion diagonale}
\hspace{20pt}

\begin{de} Les $q$-analogues respectifs de $n, n!, \binom{n}{k}$ sont les polynômes à coefficients entiers positifs:
\[  [n]_q:=1+q+q^2+ \cdots +q^{n-1}, \]
\[                   [n]!_q:=[1]_q\cdot[2]_q\cdots[n]_q,\]
\[       \binom{n}{k}_q:=\frac{[n]!_q}{[k]!_q[n-k]!_q}. \]

De plus,  nous définissons le symbole de Pochamer par:
\[ (a;q)_n:=(1-a)(1-aq)\cdots(1-aq^{n-1}) \]
\end{de}

\begin{de} Pour $\alpha=a_1a_2\cdots a_n \in \sch_{n,d}$, la statistique d'inversions diagonales notée $\mathbf{\dinv(\alpha)}$ ($r$=1) est le cardinal de l'ensemble: 
\begin{align*} \dinv(\alpha)=|\{ (i,j)|~ i<j,  &~a_i \text{ pas barré et } \A{_i}(\alpha)=\A{_j}(\alpha) \text{ ou } 
\\							      &~a_j \text{ pas barré et } \A{_i}(\alpha)=\A{_j}(\alpha)+1 \}|. 
\end{align*}

\end{de}

De façon équivalente:
\[ \dinv(\alpha)=\begin{cases} 1  & \text{ si } \A{_i}(\alpha)-\A{_j}(\alpha)=0 , a_i \text{ pas barré,}  a_j \text{ barré  et } i<j,
\\					      1  & \text{ si } 0 \geq \A{_i}(\alpha)-\A{_j}(\alpha) \geq 1 , a_i \text{ pas barré,}  a_j \text{ pas barré et } i<j,
\\					      1  & \text{ si } \A{_i}(\alpha)-\A{_j}(\alpha)=1 , a_i  \text{ barré,}  a_j \text{ pas barré et } i<j.
\end{cases} \]
 Il a été démontré dans \cite{2}\footnote{ Dans la référence $d$ est le nombre de pas vers le bas.}, aux pages 10-12, que:

\begin{align*} \s{_{n,d}}(q,1/q)& =q^{\binom{n}{2}-\binom{n-d}{2}}\sum_{\alpha \in \sch_{n,d}} q^{\A(\alpha)-\dinv(\alpha)}
\\
\\			 &= \frac{1}{[d+1]_q}\binom{n}{d}_q\binom{d+n}{n}_q =\frac{\big(q^{n+d};1/q\big)_n}{[d+1]_q(n-d)!_qd!_q(1-q)^n}, 
\end{align*}

\begin{de} Dans \cite{3}, à la page 49, et dans \cite{4}, à la page 9, la statistique d'inversion diagonale de Haiman pour $\alpha \in \dd_n^{~r}$, est défini comme:

\[ \dinv(\alpha)=\begin{cases}  r-\A_i(\alpha)+\A_j(\alpha)+1 &  \text{ si } 1 \leq \A_i(\alpha) -\A_j(\alpha) \leq r , i<j,
\\                                             r+\A_i(\alpha)-\A_j(\alpha)  & \text{ si } -r+1 \leq \A_i(\alpha) -\A_j(\alpha) \leq 0, i<j,
\\                              0& \text{ sinon. }
                        \end{cases}   \]
\end{de}

Dans \cite{3} on discute du fait que:
\[  C_n^r(q,1/q)= \sum_{\alpha \in \dd_n^{~r}} q^{\A(\alpha)-\dinv(\alpha)} \]

Notons que $C_n^1(q,1/q)=\s_{n,n}(q,1/q)$.
\\

 Soit $\alpha \in \sch_{n,d}^r$, avec $\ag(\alpha)=a_1a_2\cdots a_n$. Nous cherchons une extension de $\dinv$ tel que nous ayons:
\[ \s{_{n,d}^r}(q,t) =\sum_{\alpha \in \sch_{n,d}^r} q^{\A(\alpha)}t^{\dinv(\alpha)} ,\]
\[  \s{_{n,d}^r}(q,1/q) =q^w\sum_{\alpha \in \sch_{n,d}^r} q^{\A(\alpha)-\dinv(\alpha)}=\frac{1}{[dr+1]_q}\binom{n}{d}_q\binom{dr+n}{n} _q ,\]
où $w={\binom{n}{2}-\binom{n-d}{2}}$ si $r=1$ et $w=0$ si $d=n$.
\\ 

Dans le cas $d=1$ des expérimentations suggèrent de poser:

\[ \dinv(\alpha)=\begin{cases} 1  & \text{ si } \A{_i}(\alpha)-\A{_j}(\alpha)=0 , a_i \text{ pas barré,}  a_j \text{ barré  et } i<j,
\\					      0  & \text{ si } -1 \geq \A{_i}(\alpha)-\A{_j}(\alpha) \geq -r+1 , a_i \text{ pas barré,}  a_j \text{ barré et } i<j,
\\					      0  & \text{ si } \A{_i}(\alpha)-\A{_j}(\alpha)=0 , a_i  \text{ barré,}  a_j \text{ pas barré et } i<j.
\end{cases} \]

\begin{rmq}Dans le cas $d=1$ si $a_i$ est barré et $a_j$ n'est pas barré alors: 
\[ \A{_i}(\alpha)-\A{_j}(\alpha)=0,\]
 est la seule possibilité. Donc, la définition due $\dinv$ est équivalente à celle de Haglund pour $r=1$. Le cas où $a_i$ et $a_j$ ne sont pas barrés  ne se produit pas, donc la définition du $\dinv$ de Haiman n'est pas impliquée. Puisque $d=0$ est trivial et  $d=n$ est défini par Haiman, nous avons alors une définition qui fonctionne pour tout $r$ et tout $d$ lorsque $n=2$.
\end{rmq}

Malheureusement cette définition ne nous donne pas:
\[ \sum_{\alpha \in \sch_{n,1}^r} q^{\A(\alpha)}=\sum_{\alpha \in \sch_{n,1}^r} q^{\dinv(\alpha)}. \]
Nous n'avons donc pas non plus $\s{_{n,1}^r}(q,t)=\s{_{n,1}^r}(t,q)$.
\\

Nous verrons à la section suivante qu'il n'existe aucune extension de $\dinv$ tel que $\s{_{n,1}^r}(q,t)=\s{_{n,1}^r}(t,q)$ pour tout $r$.
\\

\section{Fonctions de stationnement 
\\sur les chemins de Schröder à pente entière.}
\hspace{20pt}

\begin{de} Une \begin{bf} $r$-fonction de stationnement  \end{bf} est une suite de longueur $n$ contenant des $k$ ou $\bar{k}$, $k \in \mathbb{N}$, possiblement distinct. Pour lequel il existe un chemin $\alpha \in  \sch_{n,d}^r$ correspondant au réordonnement croissant de la suite. 
\end{de}

\begin{note} L'ensemble des fonctions stationnement associées au chemin $\alpha$ est noté $\p(\alpha)$ et considérons  les ensembles:
 \[\p{_{n,d}^r}=\bigcup\limits_{\alpha \in \sch_{n,d}^r} \p(\alpha) \hspace{20pt} \text{ et } \hspace{20pt} \p{_{n}^r}=\bigcup_{d=0}^n \p{_{n,d}^r}. \]
\end{note}

Pour deux chemins $\alpha$ et $\beta$ distinct les ensembles $\p(\alpha)$ et $\p(\beta)$ sont clairement disjoints.
De plus, la cardinalité de $\p{_{n,d}^r}$ est donnée par la jolie formule:
\[ \frac{n!}{d!}\binom{d(r-1)+n+1}{n-d}(d(r-1)+n+1)^{d-1}. \]

Les $r$-fonctions de stationnement peuvent être vues comme des fonctions de préférences ayant $d$ places réservées. Dont l'ordre correspond à l'ordre d'arrivée. Donc pour chaque fonction de stationnement, disons $\mathcal{P}$, nous pouvons associé un unique chemin $\alpha$ et une unique permutation $\sigma$ tel que $\mathcal{P}=\sigma(\alpha)$.

\begin{tikzpicture}[line cap=round,line join=round,>=triangle 45,x=1.0cm,y=1.0cm]
\clip(-3.7530578512396677,-2.0800000000000067) rectangle (10.366942148760337,6.080000000000002);
\fill[color=xfqqff,fill=xfqqff,fill opacity=0.1] (0.0,-0.0) -- (0.0,5.0) -- (10.0,0.0) -- cycle;
\fill[color=ffcqcb,fill=ffcqcb,fill opacity=0.45] (0.0,3.0) -- (2.0,3.0) -- (2.0,2.0) -- (0.0,2.0) -- cycle;
\fill[color=ffcqcb,fill=ffcqcb,fill opacity=0.45] (0.0,2.0) -- (0.0,1.0) -- (4.0,1.0) -- (4.0,1.92) -- cycle;
\fill[color=ffcqcb,fill=ffcqcb,fill opacity=0.45] (0.0,1.0) -- (0.0,-0.0) -- (4.0,0.0) -- (4.0,1.02) -- cycle;
\draw [color=xfqqff] (0.0,-0.0)-- (0.0,5.0);
\draw [color=xfqqff] (0.0,5.0)-- (10.0,0.0);
\draw [color=xfqqff] (10.0,0.0)-- (0.0,-0.0);
\draw [line width=5.2pt,color=yqqqqq] (0.0,5.0)-- (0.0,4.0);
\draw [line width=5.2pt,color=yqqqqq] (0.0,4.0)-- (2.0,3.0);
\draw [line width=5.2pt,color=yqqqqq] (4.0,1.92)-- (4.0,1.0);
\draw [line width=5.2pt,color=yqqqqq] (4.0,1.0)-- (4.0,1.02);
\draw [line width=5.2pt,color=yqqqqq] (7.0,0.0)-- (10.0,0.0);
\draw [color=ffcqcb] (0.0,3.0)-- (2.0,3.0);
\draw [color=ffcqcb] (2.0,3.0)-- (2.0,2.0);
\draw [color=ffcqcb] (2.0,2.0)-- (0.0,2.0);
\draw [color=ffcqcb] (0.0,2.0)-- (0.0,3.0);
\draw [color=ffcqcb] (0.0,2.0)-- (0.0,1.0);
\draw [color=ffcqcb] (0.0,1.0)-- (4.0,1.0);
\draw [color=ffcqcb] (4.0,1.0)-- (4.0,1.92);
\draw [color=ffcqcb] (4.0,1.92)-- (0.0,2.0);
\draw [color=ffcqcb] (0.0,1.0)-- (0.0,-0.0);
\draw [color=ffcqcb] (0.0,-0.0)-- (4.0,0.0);
\draw [color=ffcqcb] (4.0,0.0)-- (4.0,1.02);
\draw [color=ffcqcb] (4.0,1.02)-- (0.0,1.0);
\draw (-3.6730578512396677,5.800000000000002) node[anchor=north west] {Ligne};
\draw (-3.353057851239668,4.86) node[anchor=north west] {$1$};
\draw (-3.353057851239668,3.8799999999999994) node[anchor=north west] {$2$};
\draw (-3.353057851239668,2.8999999999999986) node[anchor=north west] {$3$};
\draw (-3.353057851239668,1.8799999999999975) node[anchor=north west] {$4$};
\draw (-3.353057851239668,0.8999999999999964) node[anchor=north west] {$5$};
\draw [color=xfqqff] (0.0,4.0)-- (2.0,4.0);
\draw [color=xfqqff] (0.0,3.0)-- (4.0,3.0);
\draw [color=xfqqff] (0.0,2.0)-- (6.0,2.0);
\draw [color=xfqqff] (0.0,1.0)-- (8.0,1.0);
\draw [color=xfqqff] (1.0319999999999998,4.484)-- (1.0,0.0);
\draw [color=xfqqff] (2.0,4.0)-- (2.0,0.0);
\draw [color=xfqqff] (3.0239999999999996,3.4880000000000004)-- (3.0,0.0);
\draw [color=xfqqff] (4.0,3.0)-- (4.0,0.0);
\draw [color=xfqqff] (6.0,2.0)-- (6.0,0.0);
\draw [color=xfqqff] (7.008,1.496)-- (7.0,0.0);
\draw [color=xfqqff] (8.0,1.0)-- (8.0,0.0);
\draw [color=xfqqff] (9.008,0.49600000000000044)-- (9.0,0.0);
\draw (-2.113057851239667,5.800000000000002) node[anchor=north west] {Aire gauche};
\draw (-1.333057851239667,4.9) node[anchor=north west] {$0$};
\draw (-1.333057851239667,3.8999999999999995) node[anchor=north west] {$\bar{0}$};
\draw (-1.333057851239667,2.8999999999999986) node[anchor=north west] {$\bar{2}$};
\draw (-1.333057851239667,1.9399999999999975) node[anchor=north west] {$4$};
\draw (-1.313057851239667,0.8999999999999964) node[anchor=north west] {$4$};
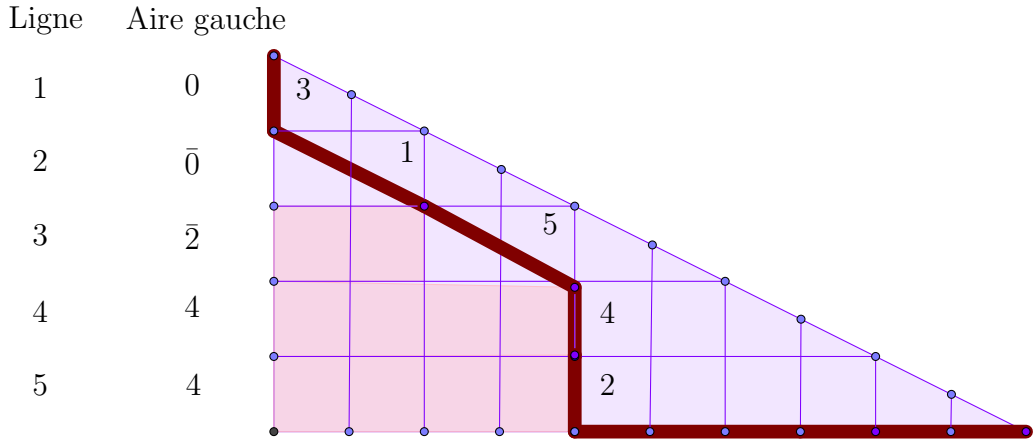
\draw (-2.5130578512396675,-0.2200000000000054) node[anchor=north west] {\footnotesize{$\text{Chemin de Schröder associé à la fonction de stationnement } \mathcal{P}=\bar{0}404\bar{2} \in Sch_{5,2}^2 $}};
\draw (-0.5130578512396675,-0.7200000000000054) node[anchor=north west] {\footnotesize{$ \text{ étiqueté par } \sigma \in S_5 \text{ tel que } \mathcal{P}=\sigma(\alpha) , \sigma=31542.$}};
\draw [line width=5.2pt,color=yqqqqq] (2.0,3.0)-- (4.0,1.92);
\draw (1.526942148760334,3.9999999999999996) node[anchor=north west] {$1$};
\draw (4.186942148760335,0.8999999999999964) node[anchor=north west] {$2$};
\draw [line width=5.2pt,color=yqqqqq] (4.0,0.0)-- (7.0,0.0);
\draw (0.14694214876033362,4.840000000000001) node[anchor=north west] {$3$};
\draw (4.186942148760335,1.8599999999999974) node[anchor=north west] {$4$};
\draw (3.4269421487603346,3.0399999999999987) node[anchor=north west] {$5$};
\draw [line width=0.4pt,color=xfqqff] (5.032,2.484)-- (5.0,0.0);
\draw [line width=5.2pt,color=yqqqqq] (4.0,1.02)-- (4.0,0.0);
\begin{scriptsize}
\draw [fill=uuuuuu] (0.0,-0.0) circle (1.5pt);
\draw [fill=xdxdff] (0.0,5.0) circle (1.5pt);
\draw [fill=xfqqff] (10.0,0.0) circle (1.5pt);
\draw [fill=xdxdff] (0.0,4.0) circle (1.5pt);
\draw [fill=xfqqff] (2.0,3.0) circle (1.5pt);
\draw [fill=xfqqff] (4.0,1.92) circle (1.5pt);
\draw [fill=qqqqff] (4.0,1.0) circle (1.5pt);
\draw [fill=xfqqff] (4.0,1.02) circle (1.5pt);
\draw [fill=xdxdff] (7.0,0.0) circle (1.5pt);
\draw [fill=xdxdff] (0.0,3.0) circle (1.5pt);
\draw [fill=xdxdff] (0.0,2.0) circle (1.5pt);
\draw [fill=xdxdff] (0.0,1.0) circle (1.5pt);
\draw [fill=xdxdff] (4.0,0.0) circle (1.5pt);
\draw [fill=xdxdff] (2.0,4.0) circle (1.5pt);
\draw [fill=xdxdff] (4.0,3.0) circle (1.5pt);
\draw [fill=xdxdff] (6.0,2.0) circle (1.5pt);
\draw [fill=xdxdff] (8.0,1.0) circle (1.5pt);
\draw [fill=xdxdff] (1.0319999999999998,4.484) circle (1.5pt);
\draw [fill=xdxdff] (1.0,0.0) circle (1.5pt);
\draw [fill=xdxdff] (2.0,0.0) circle (1.5pt);
\draw [fill=xdxdff] (3.0239999999999996,3.4880000000000004) circle (1.5pt);
\draw [fill=xdxdff] (3.0,0.0) circle (1.5pt);
\draw [fill=xdxdff] (4.0,0.0) circle (1.5pt);
\draw [fill=xdxdff] (5.032,2.484) circle (1.5pt);
\draw [fill=xdxdff] (6.0,0.0) circle (1.5pt);
\draw [fill=xdxdff] (7.008,1.496) circle (1.5pt);
\draw [fill=xfqqff] (8.0,0.0) circle (1.5pt);
\draw [fill=xdxdff] (9.008,0.49600000000000044) circle (1.5pt);
\draw [fill=xdxdff] (9.0,0.0) circle (1.5pt);
\draw [fill=xdxdff] (5.0,0.0) circle (1.5pt);
\end{scriptsize}
\end{tikzpicture}

\begin{prop} Il n'existe aucune définition de $\dinv$ sur les $r$-Schröder, rendant $\s{_{n,d}^r}(q,t)$ symétrique pour tout $r$.  
\end{prop}

\begin{pre} Par contradiction, supposons que $\s{_{2,1}^2}(q,t)=\s{_{2,1}^2}(t,q)$. 

Puisque les chemins de $\sch_{2,1}^2$ sont :
\\

\definecolor{yqqqqq}{rgb}{0.5019607843137255,0.0,0.0}
\definecolor{xfqqff}{rgb}{0.4980392156862745,0.0,1.0}
\definecolor{qqqqff}{rgb}{0.0,0.0,1.0}
\definecolor{xdxdff}{rgb}{0.49019607843137253,0.49019607843137253,1.0}
\definecolor{uuuuuu}{rgb}{0.26666666666666666,0.26666666666666666,0.26666666666666666}
\begin{tikzpicture}[line cap=round,line join=round,>=triangle 45,x=0.5cm,y=0.5cm]
\clip(-7.335030864197537,-1.113755144032923) rectangle (14.577469135802492,2.5862448559670783);
\fill[color=xfqqff,fill=xfqqff,fill opacity=0.1] (0.0,-0.0) -- (0.0,2.0125) -- (4.0,0.0125) -- cycle;
\fill[color=xfqqff,fill=xfqqff,fill opacity=0.1] (10.0,2.0) -- (10.0,0.0) -- (14.0,0.0) -- cycle;
\fill[color=xfqqff,fill=xfqqff,fill opacity=0.1] (4.914969135802479,1.9862448559670782) -- (4.914969135802479,-0.013755144032922634) -- (8.91496913580249,-0.013755144032922634) -- cycle;
\fill[color=xfqqff,fill=xfqqff,fill opacity=0.1] (-4.9975308641975404,1.998744855967078) -- (-4.9975308641975404,-0.0012551440329211516) -- (-0.997530864197532,-0.0012551440329211516) -- cycle;
\draw [color=xfqqff] (0.0,-0.0)-- (0.0,2.0125);
\draw [color=xfqqff] (0.0,2.0125)-- (4.0,0.0125);
\draw [color=xfqqff] (4.0,0.0125)-- (0.0,-0.0);
\draw [color=xfqqff] (0.0,1.018828125)-- (1.998193415637866,1.0134032921810672);
\draw [color=xfqqff] (1.996366255144038,1.0143168724279812)-- (1.9969135802469191,0.006240354938271623);
\draw [color=xfqqff] (3.0074773662551495,0.5087613168724254)-- (3.002469135802476,0.009382716049382738);
\draw (.5274691358024725,-0.4387551440329227) node[anchor=north west] {\tiny{Aire : 0}};
\draw [line width=5.2pt,color=yqqqqq] (0.0,2.0125)-- (0.0,1.018828125);
\draw [line width=5.2pt,color=yqqqqq] (0.0,1.018828125)-- (1.9963667484708045,1.0134082513426375);
\draw [color=xfqqff] (1.0024691358024733,0.0031327160493827295)-- (0.9985884773662589,1.5132057613168708);
\draw [line width=5.2pt,color=yqqqqq] (1.996366255144038,1.0143168724279812)-- (4.0,0.0125);
\draw [color=xfqqff] (10.0,2.0)-- (10.0,0.0);
\draw [color=xfqqff] (10.0,0.0)-- (14.0,0.0);
\draw [color=xfqqff] (14.0,0.0)-- (10.0,2.0);
\draw [line width=0.4pt,color=xfqqff] (10.997477366255154,1.501261316872423)-- (11.002469135802482,0.0);
\draw [line width=0.4pt,color=xfqqff] (11.992477366255155,1.0037613168724224)-- (12.002469135802482,0.0);
\draw [line width=0.4pt,color=xfqqff] (12.992477366255155,0.5037613168724224)-- (12.989969135802491,-0.0012551440329225732);
\draw [line width=0.4pt,color=xfqqff] (10.0,0.9987448559670782)-- (11.997477366255154,1.001261316872423);
\draw [color=xfqqff] (4.914969135802479,1.9862448559670782)-- (4.914969135802479,-0.013755144032922634);
\draw [color=xfqqff] (4.914969135802479,-0.013755144032922634)-- (8.91496913580249,-0.013755144032922634);
\draw [color=xfqqff] (8.91496913580249,-0.013755144032922634)-- (4.914969135802479,1.9862448559670782);
\draw [line width=0.4pt,color=xfqqff] (5.9124465020576356,1.487506172839501)-- (5.917438271604962,-0.013755144032922634);
\draw [line width=0.4pt,color=xfqqff] (6.90744650205764,0.9900061728395)-- (6.917438271604967,-0.013755144032922634);
\draw [line width=0.4pt,color=xfqqff] (7.907446502057644,0.49000617283949943)-- (7.902469135802484,-0.013755144032922577);
\draw [line width=0.4pt,color=xfqqff] (4.914969135802479,0.9849897119341566)-- (6.912446502057639,0.9875061728395005);
\draw [color=xfqqff] (-4.9975308641975404,1.998744855967078)-- (-4.9975308641975404,-0.0012551440329211516);
\draw [color=xfqqff] (-4.9975308641975404,-0.0012551440329211516)-- (-0.997530864197532,-0.0012551440329211516);
\draw [color=xfqqff] (-0.997530864197532,-0.0012551440329211516)-- (-4.9975308641975404,1.998744855967078);
\draw [line width=0.4pt,color=xfqqff] (-4.000053497942385,1.5000061728395016)-- (-3.995061728395057,-0.0012551440329211516);
\draw [line width=0.4pt,color=xfqqff] (-3.005053497942382,1.0025061728395013)-- (-2.9950617283950542,-0.0012551440329211516);
\draw [line width=0.4pt,color=xfqqff] (-2.0050534979423786,0.5025061728395008)-- (-2.010030864197532,-0.026255144032922582);
\draw [line width=0.4pt,color=xfqqff] (-4.9975308641975404,0.997489711934157)-- (-3.000053497942383,1.0000061728395018);
\draw [line width=5.2pt,color=yqqqqq] (10.0,2.0)-- (10.0,0.9987448559670782);
\draw [line width=5.2pt,color=yqqqqq] (10.0,0.9987448559670782)-- (12.002469135802482,0.0);
\draw [line width=5.2pt,color=yqqqqq] (12.002469135802482,0.0)-- (14.0,0.0);
\draw [line width=5.2pt,color=yqqqqq] (4.914969135802479,1.9862448559670782)-- (4.914969135802479,0.9849897119341566);
\draw [line width=5.2pt,color=yqqqqq] (4.914969135802479,0.9849897119341566)-- (5.9141132092408615,0.986248453105441);
\draw [line width=5.2pt,color=yqqqqq] (5.9141132092408615,0.986248453105441)-- (7.902469135802484,-0.013755144032922577);
\draw [line width=5.2pt,color=yqqqqq] (7.902469135802484,-0.013755144032922577)-- (8.91496913580249,-0.013755144032922634);
\draw [line width=5.2pt,color=yqqqqq] (-4.9975308641975404,1.998744855967078)-- (-3.000053497942383,1.0000061728395018);
\draw [line width=5.2pt,color=yqqqqq] (-3.000053497942383,1.0000061728395018)-- (-2.9950617283950542,-0.0012551440329211516);
\draw [line width=5.2pt,color=yqqqqq] (-2.9950617283950542,-0.0012551440329211516)-- (-0.997530864197532,-0.0012551440329211516);
\draw (-4.310030864197535,-0.4387551440329227) node[anchor=north west] {\tiny{Aire : 0}};
\draw (5.58996913580248,-0.4387551440329227) node[anchor=north west] {\tiny{Aire : 1}};
\draw (10.714969135802487,-0.4387551440329227) node[anchor=north west] {\tiny{Aire : 2}};
\begin{scriptsize}
\draw [fill=uuuuuu] (0.0,-0.0) circle (1.5pt);
\draw [fill=xdxdff] (0.0,2.0125) circle (1.5pt);
\draw [fill=qqqqff] (4.0,0.0125) circle (1.5pt);
\draw [fill=xdxdff] (0.0,1.018828125) circle (1.5pt);
\draw [fill=xdxdff] (1.998193415637866,1.0134032921810672) circle (1.5pt);
\draw [fill=xdxdff] (1.996366255144038,1.0143168724279812) circle (1.5pt);
\draw [fill=xdxdff] (1.9969135802469191,0.006240354938271623) circle (1.5pt);
\draw [fill=xdxdff] (1.0024691358024733,0.0031327160493827295) circle (1.5pt);
\draw [fill=xdxdff] (3.0074773662551495,0.5087613168724254) circle (1.5pt);
\draw [fill=xdxdff] (3.002469135802476,0.009382716049382738) circle (1.5pt);
\draw [fill=uuuuuu] (1.9963667484708045,1.0134082513426375) circle (1.5pt);
\draw [fill=xdxdff] (0.9985884773662589,1.5132057613168708) circle (1.5pt);
\draw [fill=qqqqff] (10.0,2.0) circle (1.5pt);
\draw [fill=qqqqff] (10.0,0.0) circle (1.5pt);
\draw [fill=qqqqff] (14.0,0.0) circle (1.5pt);
\draw [fill=xdxdff] (10.997477366255154,1.501261316872423) circle (1.5pt);
\draw [fill=xdxdff] (11.002469135802482,0.0) circle (1.5pt);
\draw [fill=xdxdff] (11.992477366255155,1.0037613168724224) circle (1.5pt);
\draw [fill=xdxdff] (12.002469135802482,0.0) circle (1.5pt);
\draw [fill=xdxdff] (12.992477366255155,0.5037613168724224) circle (1.5pt);
\draw [fill=qqqqff] (12.989969135802491,-0.0012551440329225732) circle (1.5pt);
\draw [fill=xdxdff] (10.0,0.9987448559670782) circle (1.5pt);
\draw [fill=xdxdff] (11.997477366255154,1.001261316872423) circle (1.5pt);
\draw [fill=qqqqff] (4.914969135802479,1.9862448559670782) circle (1.5pt);
\draw [fill=qqqqff] (4.914969135802479,-0.013755144032922634) circle (1.5pt);
\draw [fill=qqqqff] (8.91496913580249,-0.013755144032922634) circle (1.5pt);
\draw [fill=xdxdff] (5.9124465020576356,1.487506172839501) circle (1.5pt);
\draw [fill=xdxdff] (5.917438271604962,-0.013755144032922634) circle (1.5pt);
\draw [fill=xdxdff] (6.90744650205764,0.9900061728395) circle (1.5pt);
\draw [fill=xdxdff] (6.917438271604967,-0.013755144032922634) circle (1.5pt);
\draw [fill=xdxdff] (7.907446502057644,0.49000617283949943) circle (1.5pt);
\draw [fill=qqqqff] (7.902469135802484,-0.013755144032922577) circle (1.5pt);
\draw [fill=xdxdff] (4.914969135802479,0.9849897119341566) circle (1.5pt);
\draw [fill=xdxdff] (6.912446502057639,0.9875061728395005) circle (1.5pt);
\draw [fill=qqqqff] (-4.9975308641975404,1.998744855967078) circle (1.5pt);
\draw [fill=qqqqff] (-4.9975308641975404,-0.0012551440329211516) circle (1.5pt);
\draw [fill=qqqqff] (-0.997530864197532,-0.0012551440329211516) circle (1.5pt);
\draw [fill=xdxdff] (-4.000053497942385,1.5000061728395016) circle (1.5pt);
\draw [fill=xdxdff] (-3.995061728395057,-0.0012551440329211516) circle (1.5pt);
\draw [fill=xdxdff] (-3.005053497942382,1.0025061728395013) circle (1.5pt);
\draw [fill=xdxdff] (-2.9950617283950542,-0.0012551440329211516) circle (1.5pt);
\draw [fill=xdxdff] (-2.0050534979423786,0.5025061728395008) circle (1.5pt);
\draw [fill=qqqqff] (-2.010030864197532,-0.026255144032922582) circle (1.5pt);
\draw [fill=xdxdff] (-4.9975308641975404,0.997489711934157) circle (1.5pt);
\draw [fill=xdxdff] (-3.000053497942383,1.0000061728395018) circle (1.5pt);
\draw [fill=uuuuuu] (5.9141132092408615,0.986248453105441) circle (1.5pt);
\end{scriptsize}
\end{tikzpicture}
\\

Nous aurions alors: 
\[ 2+q+q^2=\sum_{\alpha \in \sch_{2,1}^2}q^{\A(\alpha)}=\s{_{2,1}^2}(q,1)=\s{_{2,1}^2}(1,q)=\sum_{\alpha \in \sch_{2,1}^2}q^{\dinv(\alpha)}.\]

Donc en particulier nous avons l'égalité de multi ensembles suivante:
 \[\{\A(\alpha) | \alpha \in \sch{_{2,1}^2 \}}=\{\dinv(\alpha) | \alpha \in \sch{_{2,1}^2 \}}=\{0,0,1,2\}.\] 

Posons: \[\A(\alpha_1)=0,~\A(\alpha_2)=0,~\A(\alpha_3)=1,~\A(\alpha_4)=2.\]
\\

Puisque nous pouvons considérer les fonctions de stationnement comme étant la permutation de multi ensembles ayant la propriété que:
\[ |\{ \mathcal{P}(i)=k : i \in [0, \cdots, n ] \}| \leq n-k.\]
 où $\mathcal{P} \in \p$ de longueur $n$ et $\mathcal{P}(i)$ donne la valeur à la position $i+1$. Puisque les chemins à pente entière ont toujours cette propriété nous pouvons alors considérer $0012$ dans $\sch_{4,4}^1$ et nous avons alors que:

\[(\dinv(\alpha_1),\dinv(\alpha_2),\dinv(\alpha_3),\dinv(\alpha_4)) \in \p(0012).\]
\\
 
Nous trouvons la contradiction en vérifiant les sur 12 possibilités $\mathcal{P} \in \p(0012)$ que:
\[q^w(\sum_{i=1}^4 q^{\A(\alpha_i)-\mathcal{P}(i)})\not=  1+q+q^2+q^3=[4]_q=\frac{1}{[2+1]_q}\binom{2}{1}_q\binom{2+2}{2} _q , w \in \mathbb{Z} \]

Cette vérification est faite à \hyperref[Annexe 1]{l'annexe 1}.
\end{pre} 

\newpage

\section{Fonctions de stationnement sur les chemins 
\\ sans diagonale incluse dans un pentagone rectangle.}

 \begin{de} Un chemin de Schröder ayant une fraction unitaire comme pente, nommée \begin{bf} $r$-FSchröder \end{bf} est une suite de pas débutant en $(0,nr)$ et se terminant en $(n,0)$ sans jamais passer au-dessus de la droite passant par $(0,nr)$ et $(n,0)$. De plus, les seuls pas autorisé sont $(0,-1)$, dit vers le bas (certains auteurs disent aussi vers sud), $(1,0)$ dits vers la droite (vers l'est) ou $(1,-r)$ un pas diagonal. \end{de}

Dans un chemin de Schröder, la forme entre deux diagonales non consécutive est un pentagone rectangle possiblement dégénéré. Afin d'établir une formule comptant les chemins de stationnement dans les $r$-FSchröder, nous allons dans cette section établir une formule comptant le nombre de chemins de stationnements sans diagonale sont dans un pentagone. Nous cherchons donc combien de permutations sont possibles sur la donnée d'aire gauche des chemins contenus dans un pentagone et composés uniquement de pas vers le bas et de pas vers la droite.
\\

Considérons la forme définie par:
 \[ P=\{(0,0),(0,a),(b,0),(b,q),(p,a)\} \]
 \[\text{ tel que: }a,b,q,r \in \mathbb{N} ,p \in \mathbb{R}^+, \frac{1}{r}=\frac{q-a}{p-b} \text{, si $p-b\not=0$, $q-a\not=0$ et } a \text{ ou } b \text{ non nul.}\]
  Par abus nous nommerons cette forme un \label{pentagone}{pentagone}, mais plusieurs cas sont possibles:
\\Cas 1 un segment vertical ou un segment horizontal.
\\Cas 2 un rectangle.
\\Cas 3 un triangle, un trapèze ou un pentagone dans ce cas notons $1/r$ la pente de la droite passant par $(p,a),(b,q)$.

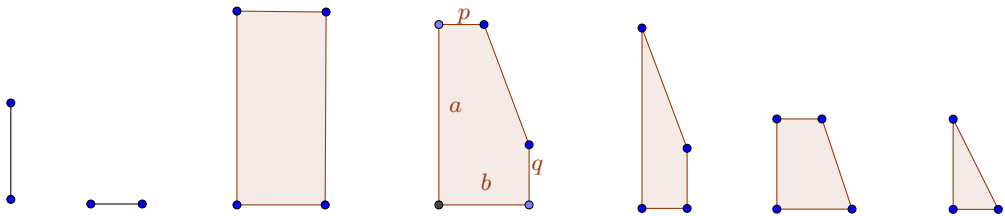
\definecolor{zzttqq}{rgb}{0.6,0.2,0.0}
\definecolor{qqqqff}{rgb}{0.0,0.0,1.0}
\definecolor{xdxdff}{rgb}{0.49019607843137253,0.49019607843137253,1.0}
\definecolor{uuuuuu}{rgb}{0.26666666666666666,0.26666666666666666,0.26666666666666666}
\begin{tikzpicture}[line cap=round,line join=round,>=triangle 45,x=0.4cm,y=0.4cm]
\clip(-15.29580249979817,-1.002558366503328) rectangle (19.86857124210292,7.051446461602022);
\fill[color=zzttqq,fill=zzttqq,fill opacity=0.1] (0.0,-0.0) -- (0.0,6.0) -- (1.5,6.0) -- (3.0,2.0) -- (2.99875678058222,0.0) -- cycle;
\fill[color=zzttqq,fill=zzttqq,fill opacity=0.1] (6.751148858292799,5.883549540333181) -- (6.751148858292799,-0.11645045966681922) -- (8.251148858292797,-0.11645045966681922) -- (8.251148858292797,1.8835495403331808) -- cycle;
\fill[color=zzttqq,fill=zzttqq,fill opacity=0.1] (11.232919952628997,2.8593738055713827) -- (11.232919952628997,-0.14062619442861485) -- (13.732919952628999,-0.14062619442861485) -- (12.732919952628999,2.8593738055713827) -- cycle;
\fill[color=zzttqq,fill=zzttqq,fill opacity=0.1] (17.0930420743381,-0.1476228844121527) -- (17.0930420743381,2.8523771155878483) -- (18.5930420743381,-0.1476228844121527) -- cycle;
\fill[color=zzttqq,fill=zzttqq,fill opacity=0.1] (-6.710929327361846,6.446093557929162) -- (-6.710929327361847,-1.9910891807746575E-4) -- (-3.7780189196967697,-1.9910891807717324E-4) -- (-3.747467769616925,6.415542407849317) -- cycle;
\draw [color=zzttqq] (0.0,-0.0)-- (0.0,6.0);
\draw [color=zzttqq] (0.0,6.0)-- (1.5,6.0);
\draw [color=zzttqq] (1.5,6.0)-- (3.0,2.0);
\draw [color=zzttqq] (3.0,2.0)-- (2.99875678058222,0.0);
\draw [color=zzttqq] (2.99875678058222,0.0)-- (0.0,-0.0);
\draw [color=zzttqq] (6.751148858292799,5.883549540333181)-- (6.751148858292799,-0.11645045966681922);
\draw [color=zzttqq] (6.751148858292799,-0.11645045966681922)-- (8.251148858292797,-0.11645045966681922);
\draw [color=zzttqq] (8.251148858292797,-0.11645045966681922)-- (8.251148858292797,1.8835495403331808);
\draw [color=zzttqq] (8.251148858292797,1.8835495403331808)-- (6.751148858292799,5.883549540333181);
\draw [color=zzttqq] (11.232919952628997,2.8593738055713827)-- (11.232919952628997,-0.14062619442861485);
\draw [color=zzttqq] (11.232919952628997,-0.14062619442861485)-- (13.732919952628999,-0.14062619442861485);
\draw [color=zzttqq] (13.732919952628999,-0.14062619442861485)-- (12.732919952628999,2.8593738055713827);
\draw [color=zzttqq] (12.732919952628999,2.8593738055713827)-- (11.232919952628997,2.8593738055713827);
\draw [color=zzttqq] (17.0930420743381,-0.1476228844121527)-- (17.0930420743381,2.8523771155878483);
\draw [color=zzttqq] (17.0930420743381,2.8523771155878483)-- (18.5930420743381,-0.1476228844121527);
\draw [color=zzttqq] (18.5930420743381,-0.1476228844121527)-- (17.0930420743381,-0.1476228844121527);
\draw (-14.22651224700361,3.3909785499446894)-- (-14.22651224700361,0.18310779156099005);
\draw (-11.568562190057133,0.030352041161766413)-- (-9.857697785585836,0.030352041161766413);
\draw [color=zzttqq] (-6.710929327361846,6.446093557929162)-- (-6.710929327361847,-1.9910891807746575E-4);
\draw [color=zzttqq] (-6.710929327361847,-1.9910891807746575E-4)-- (-3.7780189196967697,-1.9910891807717324E-4);
\draw [color=zzttqq] (-3.7780189196967697,-1.9910891807717324E-4)-- (-3.747467769616925,6.415542407849317);
\draw [color=zzttqq] (-3.747467769616925,6.415542407849317)-- (-6.710929327361846,6.446093557929162);
\begin{scriptsize}
\draw [fill=uuuuuu] (0.0,-0.0) circle (1.5pt);
\draw [fill=xdxdff] (0.0,6.0) circle (1.5pt);
\draw [fill=qqqqff] (1.5,6.0) circle (1.5pt);
\draw [fill=qqqqff] (3.0,2.0) circle (1.5pt);
\draw [fill=xdxdff] (2.99875678058222,0.0) circle (1.5pt);
\draw[color=zzttqq] (0.5602443916411581,3.2687739496253103) node {$a$};
\draw[color=zzttqq] (0.8352047423597592,6.273968256172579) node {$p$};
\draw[color=zzttqq] (3.276620296910899,1.2523980443555571) node {$q$};
\draw[color=zzttqq] (1.5684323442760286,0.7635796430780412) node {$b$};
\draw [fill=qqqqff] (6.751148858292799,5.883549540333181) circle (1.5pt);
\draw [fill=qqqqff] (6.751148858292799,-0.11645045966681922) circle (1.5pt);
\draw [fill=qqqqff] (8.251148858292797,-0.11645045966681922) circle (1.5pt);
\draw [fill=qqqqff] (8.251148858292797,1.8835495403331808) circle (1.5pt);
\draw [fill=qqqqff] (11.232919952628997,2.8593738055713827) circle (1.5pt);
\draw [fill=qqqqff] (11.232919952628997,-0.14062619442861485) circle (1.5pt);
\draw [fill=qqqqff] (13.732919952628999,-0.14062619442861485) circle (1.5pt);
\draw [fill=qqqqff] (12.732919952628999,2.8593738055713827) circle (1.5pt);
\draw [fill=qqqqff] (17.0930420743381,2.8523771155878483) circle (1.5pt);
\draw [fill=qqqqff] (17.0930420743381,-0.1476228844121527) circle (1.5pt);
\draw [fill=qqqqff] (18.5930420743381,-0.1476228844121527) circle (1.5pt);
\draw [fill=qqqqff] (-14.22651224700361,3.3909785499446894) circle (1.5pt);
\draw [fill=qqqqff] (-14.22651224700361,0.18310779156099005) circle (1.5pt);
\draw [fill=qqqqff] (-11.568562190057133,0.030352041161766413) circle (1.5pt);
\draw [fill=qqqqff] (-9.857697785585836,0.030352041161766413) circle (1.5pt);
\draw [fill=qqqqff] (-6.710929327361846,6.446093557929162) circle (1.5pt);
\draw [fill=qqqqff] (-6.710929327361847,-1.9910891807746575E-4) circle (1.5pt);
\draw [fill=qqqqff] (-3.7780189196967697,-1.9910891807717324E-4) circle (1.5pt);
\draw [fill=qqqqff] (-3.747467769616925,6.415542407849317) circle (1.5pt);
\end{scriptsize}
\end{tikzpicture}
\\ Remarquons d'abord que dans le cas 2 l'aire gauche est comprise entre $0$ et $B$, il y a donc $B+1$ possibilités. Comme pour tout chemin $\alpha$ contenu dans ce rectangle $\ag(\alpha)$ est une suite unique de $a$ nombre, parmi $\{0,\cdots, b\}$, en ordre croissant. Nous avons $(b+1)^a$ arrangements, donc $(b+1)^a$ permutations sur l'ensemble des chemins. 
\\

De plus nous pouvons également remarquer que le cas 1 est inclue dans ce cas, car si $a=0$ nous avons $(b+1)^0=1$ et si $b=0$ nous avons $(0+1)^a=1$.
\\

Dans le cas 3, si $a \leq r$ et $q=0$ le seul pas possible dans la dernière colonne est le pas allant de $(b-1,0)$ vers $(b,0)$ puisque la pente est de $1/r \leq 1$. Donc le nombre de chemins possible est contenu dans le rectangle $(0,a), (0,0), (b-1,0), (b-1,a)$ et par ce qui précède le nombre de permutations de l'aire de ses chemins est donné par $(b-1+1)^a=b^a$.
\\

Si $a-q \leq r$ et $q\not=0$, alors pour les mêmes raisons que précédemment  les chemins sont contenus dans l'équerre $(0,a), (0,0), (b,0), (b,q), (b-1,q), (b-1,a)$. Donc tout chemin doit avoir un pas vers la droite débutant par $(b-1,i)$, où $ i \in \{0,Q\}$.
Comme l'ensemble des valeurs possible pour l'aire gauche d'une ligne avant un pas vers la droite est disjoint de l'unique valeur possible pour l'aire gauche d'une ligne après un pas vers la droite, nous pouvons faire le produit de combien de fonctions stationnements sont possible avant le pas droit avec la binomiale de comment permuter les deux ensembles. 
\\

Si $a-q > r$ et $p$ n'est pas entier la première ligne ou il est possible d'avoir un pas vers la droite dans le triangle $(p,q),(p,a),(b,q)$ est donné par la formule $a-q+(1-b+\lfloor p \rfloor)r$. Donc si le pentagone n'est pas dégénéré en un rectangle tout chemin doit avoir un pas vers la droite débutant par $(\lfloor p \rfloor,i)$, où $ i \in \{a-q+(1-b+\lfloor p \rfloor)r,\cdots ,a\}$. Clairement cette formule fonctionne également si $p$ est un entier.
\\

Comme l'ensemble des valeurs possible pour l'aire gauche d'une ligne avant un pas vers la droite est disjoint de l'ensemble des valeurs possible pour l'aire gauche d'une ligne après un pas vers la droite, nous pouvons faire le produit de combien de fonctions de stationnement sont possible avant le pas droit, combien sont possible après et comment permuter les deux ensembles. 
\\

Dans tous les cas, la partie de droite est de largeur $b-\lfloor p\rfloor-1$ qui est strictement plus petite que $b$ et de hauteur $q+(b-\lfloor p \rfloor-1)r$ qui est strictement plus petite que $a$, car autrement nous aurions:
\[ (b- p)r  > (b- p- 1)r  \geq (b-\lfloor p \rfloor-1)r \geq a-q= (b- p)r.\]
\\

Donc le processus prend au plus le maximum de $a-q$ ou $b$ étapes.

\newpage

Nous obtenons alors la formule récursive:

\[ \PP(a,b,p,r,q)=\begin{cases}  (b+1)^a  & \text{si $P$ définie un rectangle},
\\
\\					           \sum\limits_{w=0}^q b^{a-q+w} \binom{a}{q-w}&  \text{ si }~ a-q \leq r,
\\
\\                                                        \sum\limits_{w=0}^q (\lfloor p\rfloor+1)^{a-w}(b-\lfloor p \rfloor)^w \binom{a}{w} 
\\                                                         +\sum\limits_{t=q+1}^{q + (b-\lfloor p \rfloor-1)r} (\lfloor p\rfloor+1)^{a-t}\binom{a}{t}\PP(\xi) & \text{ sinon, }		   
\end{cases} \]

  où $ \xi= (t,b-\lfloor p \rfloor-1, b- \lfloor p \rfloor - \frac{t-q}{r}-1 ,r,q).$	Il importe de remarquer que si $a=0$ ou $b=0$ cette formule est bien définie, mais qu'elle ne l'est pas si $a=b=0$.
  
  \begin{prop} $ \PP$ compte le nombre de fonctions de stationnement dans un pentagone possiblement dégénéré ayant une fraction unitaire comme pente. 
  \end{prop}
  \begin{pre} Par construction $\F$
  \end{pre}
  
  De façon analogue on trouve pour un pentagone à pente entière:
  
  \[ \PPE(a,b,p,r,q)=\begin{cases}  (b+1)^a  & \text{si $P$ définie un rectangle},
\\
\\					           \sum\limits_{w=0}^q b^{a-q+w} \binom{a}{q-w}&  \text{ si }~ a-q \leq 1,
\\
\\                                                        \sum\limits_{w=0}^q ( p+1)^{a-w}(b- p)^w \binom{a}{w} 
\\                                                         +\sum\limits_{t=q+1}^{a-1} (p+1)^{a-t}\binom{a}{t}\PPE(\xi) & \text{ sinon, }		   
\end{cases} \]
  
  où $ \xi= (t,b-p-1, (a-t)r-1,r,q).$ Il importe de remarquer que si $a=0$ ou $b=0$ cette formule est bien définie, mais qu'elle ne l'est pas si $a=b=0$.
  \\
\newpage   

  \section{Fonctions de stationnement sur les $r$-FSchröder}

Soit $r \in \mathbb{N}$ telle que $1/r$ soit la pente du triangle rectangle contenant le chemin et $n$ la largeur de la base. Nous avons donc que $nr$ est la hauteur. Et soit $d$ le nombre de pas qui descendent. Il est facile de remarquer que le nombre de chemins est égal à ceux contenus dans $\sch_{n,d}^r$.
 \\

Les fonctions de stationnement sur les $r$-FSchröder (chemins de Schröder ayant une fraction unitaire comme pente), peuvent alors être vues comme des fonctions de préférences ayant $d/r$ sections réservées. Dont l'ordre correspond à l'ordre d'arrivée. Donc pour chaque fonction de stationnement, disons $\mathcal{P}$, nous pouvons associé un unique chemin $\alpha$ et une unique permutation $\sigma$ tel que $\mathcal{P}=\sigma(\alpha)$.
\\

\definecolor{ffcqcb}{rgb}{1.0,0.7529411764705882,0.796078431372549}
\definecolor{uuuuuu}{rgb}{0.26666666666666666,0.26666666666666666,0.26666666666666666}
\definecolor{yqqqqq}{rgb}{0.5019607843137255,0.0,0.0}
\definecolor{xfqqff}{rgb}{0.4980392156862745,0.0,1.0}
\begin{tikzpicture}[line cap=round,line join=round,>=triangle 45,x=.8cm,y=.8cm]
\clip(-6.9836320000000045,-1.4771592000000073) rectangle (10.632368000000006,7.0368408);
\fill[color=xfqqff,fill=xfqqff,fill opacity=0.1] (0.0,6.0) -- (0.0,-0.0) -- (3.0,0.0) -- cycle;
\fill[color=ffcqcb,fill=ffcqcb,fill opacity=1.0] (0.0,2.0) -- (2.0,2.0) -- (2.0,0.0) -- (0.0,-0.0) -- cycle;
\draw [color=xfqqff] (0.0,6.0)-- (0.0,-0.0);
\draw [color=xfqqff] (0.0,-0.0)-- (3.0,0.0);
\draw [color=xfqqff] (3.0,0.0)-- (0.0,6.0);
\draw [color=xfqqff] (0.0,5.0)-- (0.5,5.0);
\draw [color=xfqqff] (0.0,4.0)-- (1.0,4.0);
\draw [color=xfqqff] (0.0,3.0)-- (1.5,3.0);
\draw [color=xfqqff] (0.0,2.0)-- (2.0,2.0);
\draw [color=xfqqff] (0.0,1.0)-- (2.5,1.0);
\draw [color=xfqqff] (1.0,0.0)-- (1.0,4.0);
\draw [color=xfqqff] (2.0,0.0)-- (2.0,2.0);
\draw [line width=5.2pt,color=yqqqqq] (0.0,6.0)-- (0.0,4.0);
\draw [line width=5.2pt,color=yqqqqq] (0.0,4.0)-- (1.0,2.0);
\draw [line width=5.2pt,color=yqqqqq] (1.0,2.0)-- (2.0,2.0);
\draw [line width=5.2pt,color=yqqqqq] (2.0,2.0)-- (3.0,0.0);
\draw [color=ffcqcb] (0.0,2.0)-- (2.0,2.0);
\draw [color=ffcqcb] (2.0,2.0)-- (2.0,0.0);
\draw [color=ffcqcb] (2.0,0.0)-- (0.0,-0.0);
\draw [color=ffcqcb] (0.0,-0.0)-- (0.0,2.0);
\draw (-5.509632000000003,6.552840799999999) node[anchor=north west] {Rangée};
\draw (-3.199632000000001,6.552840799999999) node[anchor=north west] {Aire gauche};
\draw (-5.157632000000003,5.980840799999999) node[anchor=north west] {$1$};
\draw (-5.157632000000003,4.946840799999999) node[anchor=north west] {$2$};
\draw (-5.1796320000000025,3.9348407999999973) node[anchor=north west] {$3$};
\draw (-5.157632000000003,2.9228407999999964) node[anchor=north west] {$4$};
\draw (-5.157632000000003,1.9548407999999957) node[anchor=north west] {$5$};
\draw (-5.157632000000003,0.9428407999999947) node[anchor=north west] {$6$};
\draw (-2.6496320000000004,6.0028407999999995) node[anchor=north west] {$0$};
\draw (-2.6496320000000004,4.990840799999998) node[anchor=north west] {$0$};
\draw (-2.6496320000000004,3.9788407999999973) node[anchor=north west] {$\bar{0}$};
\draw (-2.6496320000000004,2.9668407999999964) node[anchor=north west] {$\bar{0}$};
\draw (-2.6496320000000004,1.9768407999999957) node[anchor=north west] {$\bar{2}$};
\draw (-2.6496320000000004,0.9868407999999947) node[anchor=north west] {$\bar{2}$};
\draw (0.34236800000000206,5.936840799999999) node[anchor=north west] {$6$};
\draw (0.3643680000000021,4.902840799999998) node[anchor=north west] {$2$};
\draw (0.4963680000000022,3.9568407999999975) node[anchor=north west] {$4$};
\draw (1.0463680000000026,3.0108407999999964) node[anchor=north west] {$3$};
\draw (2.498368000000004,1.9328407999999957) node[anchor=north west] {$5$};
\draw (3.0483680000000044,0.9428407999999947) node[anchor=north west] {$1$};
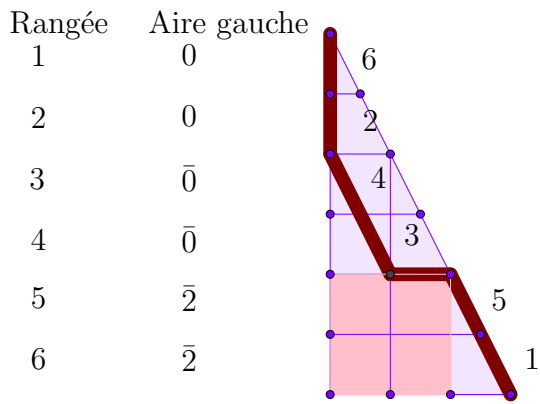
\draw (-6.565632000000004,-0.3331592000000064) node[anchor=north west] {$\text{Chemin de Schröder associé à la fonction de stationnement } \mathcal{P}=\bar{2}0\bar{0}\bar{0}\bar{2}0$};
\begin{scriptsize}
\draw [fill=xfqqff] (0.0,6.0) circle (1.5pt);
\draw [fill=xfqqff] (0.0,-0.0) circle (1.5pt);
\draw [fill=xfqqff] (3.0,0.0) circle (1.5pt);
\draw [fill=xfqqff] (0.0,5.0) circle (1.5pt);
\draw [fill=xfqqff] (0.5,5.0) circle (1.5pt);
\draw [fill=xfqqff] (0.0,4.0) circle (1.5pt);
\draw [fill=xfqqff] (1.0,4.0) circle (1.5pt);
\draw [fill=xfqqff] (0.0,3.0) circle (1.5pt);
\draw [fill=xfqqff] (1.5,3.0) circle (1.5pt);
\draw [fill=xfqqff] (0.0,2.0) circle (1.5pt);
\draw [fill=xfqqff] (2.0,2.0) circle (1.5pt);
\draw [fill=xfqqff] (0.0,1.0) circle (1.5pt);
\draw [fill=xfqqff] (2.5,1.0) circle (1.5pt);
\draw [fill=xfqqff] (1.0,0.0) circle (1.5pt);
\draw [fill=xfqqff] (2.0,0.0) circle (1.5pt);
\draw [fill=uuuuuu] (1.0,2.0) circle (1.5pt);
\end{scriptsize}
\end{tikzpicture}

Remarquons d'abord qu'il y a toujours un multiple de $r$ lignes diagonales consécutives dans les $r$-FSchröder. De plus, la valeur de l'aire gauche d'une diagonale (d'une valeur barrée) est toujours présente exactement $r$ fois. En effet, le chemin ne peut changer de direction que lorsqu'il est sur une coordonnée à valeur entière et la pente étant de $1/r$ nous devons donc avoir $r$ pas diagonales pour atteindre une autre coordonnée à valeur entière. Donc dans ce type de chemins si $d$ n'est pas un multiple de $r$ alors il y a aucun chemin. De plus, si nous avons $2r$ diagonales consécutives les $r$ premières ne peuvent pas avoir la même valeur que les autres, car $r$ diagonales équivalentes à $r$ pas qui descendent et un pas vers la droite.

\newpage

Ainsi si nous connaissons le nombre de permutations de l'aire gauche qui sont possible pour  les parties non diagonales il suffit de le multiplier par le  multinôme: 

\[ \binom{nr}{d,\underbrace{r,\cdots,r}_\textrm{n-d/r fois}} \]
\\

Remarquons que les parties non diagonales sont aux nombres d'au plus $d+d/r$, car nous devons ajouter $d$ pas vers le bas et $d/r$ pas vers la droite. Et que celles-ci sont toute des pentagones possiblement dégénérés vue à la section 4. Et sont partagées parmi $(nr-d)/r=n-d/r$ groupes de $r$ diagonales consécutives. De plus, pour un partage en $i$ pentagones ayant une diagonale entre eux, nous avons une façon de partager les $n-d/r$ groupes de diagonales de façon à commencer et terminer par des diagonales, donc en $i+1$ parties (respectivement, commencer et terminer par des sous-chemins sans diagonales donc en $i-1$ parties) et deux façons de partager les $n-d/r$ groupes de diagonales de façon à avoir des diagonales au début et un pentagone à la fin ou un pentagone au début et  des diagonales à la fin, donc en $i$ \label{id}{parties}. 

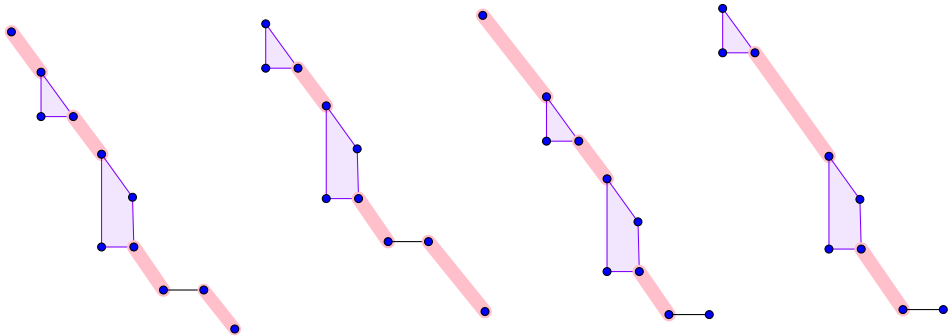
\definecolor{xfqqff}{rgb}{0.4980392156862745,0.0,1.0}
\definecolor{ffcqcb}{rgb}{1.0,0.7529411764705882,0.796078431372549}
\definecolor{qqqqff}{rgb}{0.0,0.0,1.0}
\begin{tikzpicture}[line cap=round,line join=round,>=triangle 45,x=1.0cm,y=1.0cm]
\clip(-8.115162334435636,3.3969048370880137) rectangle (6.935126604159587,9.09934995396059);
\fill[color=xfqqff,fill=xfqqff,fill opacity=0.1] (-6.462675294049269,7.405765819270167) -- (-6.033413039423397,6.815530219159594) -- (-6.462675294049269,6.815530219159594) -- cycle;
\fill[color=xfqqff,fill=xfqqff,fill opacity=0.1] (-5.657808566625761,6.314724255429412) -- (-5.246432239275968,5.742374582594918) -- (-5.2285463119998905,5.080595273380034) -- (-5.657808566625761,5.080595273380034) -- cycle;
\fill[color=xfqqff,fill=xfqqff,fill opacity=0.1] (-3.4757254389442522,8.049659201208971) -- (-3.046463184318381,7.4594236010984) -- (-3.4757254389442522,7.4594236010984) -- cycle;
\fill[color=xfqqff,fill=xfqqff,fill opacity=0.1] (-2.6708587115207436,6.958617637368217) -- (-2.2594823841709473,6.386267964533723) -- (-2.2415964568948703,5.724488655318839) -- (-2.6708587115207436,5.724488655318839) -- cycle;
\fill[color=xfqqff,fill=xfqqff,fill opacity=0.1] (0.25702791899225935,7.078827808345808) -- (0.686290173618133,6.488592208235236) -- (0.25702791899225935,6.488592208235236) -- cycle;
\fill[color=xfqqff,fill=xfqqff,fill opacity=0.1] (1.0618946464157715,5.9877862445050525) -- (1.4732709737655632,5.415436571670558) -- (1.491156901041641,4.753657262455674) -- (1.0618946464157715,4.753657262455674) -- cycle;
\fill[color=xfqqff,fill=xfqqff,fill opacity=0.1] (2.599087002017442,8.254552666897432) -- (3.028349256643314,7.664317066786861) -- (2.599087002017442,7.664317066786861) -- cycle;
\fill[color=xfqqff,fill=xfqqff,fill opacity=0.1] (4.012075256827602,6.287100666528859) -- (4.4234515841774,5.714750993694366) -- (4.441337511453479,5.052971684479483) -- (4.012075256827602,5.052971684479483) -- cycle;
\draw [line width=5.2pt,color=ffcqcb] (-6.856165694122984,7.942343637552505)-- (-6.462675294049269,7.405765819270167);
\draw [color=xfqqff] (-6.462675294049269,7.405765819270167)-- (-6.033413039423397,6.815530219159594);
\draw [color=xfqqff] (-6.033413039423397,6.815530219159594)-- (-6.462675294049269,6.815530219159594);
\draw [color=xfqqff] (-6.462675294049269,6.815530219159594)-- (-6.462675294049269,7.405765819270167);
\draw [line width=5.2pt,color=ffcqcb] (-6.033413039423397,6.815530219159594)-- (-5.657808566625761,6.314724255429412);
\draw [color=xfqqff] (-5.657808566625761,6.314724255429412)-- (-5.246432239275968,5.742374582594918);
\draw [color=xfqqff] (-5.246432239275968,5.742374582594918)-- (-5.2285463119998905,5.080595273380034);
\draw [color=xfqqff] (-5.2285463119998905,5.080595273380034)-- (-5.657808566625761,5.080595273380034);
\draw [color=xfqqff] (-5.657808566625761,5.080595273380034)-- (-5.657808566625761,6.314724255429412);
\draw [line width=5.2pt,color=ffcqcb] (-5.2285463119998905,5.080595273380034)-- (-4.835055911926176,4.50824560054554);
\draw (-4.835055911926176,4.50824560054554)-- (-4.298478093643838,4.50824560054554);
\draw [line width=4.800000000000001pt,color=ffcqcb] (-4.298478093643838,4.50824560054554)-- (-3.8871017662940446,3.989553709539279);
\draw [color=xfqqff] (-3.4757254389442522,8.049659201208971)-- (-3.046463184318381,7.4594236010984);
\draw [color=xfqqff] (-3.046463184318381,7.4594236010984)-- (-3.4757254389442522,7.4594236010984);
\draw [color=xfqqff] (-3.4757254389442522,7.4594236010984)-- (-3.4757254389442522,8.049659201208971);
\draw [line width=5.2pt,color=ffcqcb] (-3.046463184318381,7.4594236010984)-- (-2.6708587115207436,6.958617637368217);
\draw [color=xfqqff] (-2.6708587115207436,6.958617637368217)-- (-2.2594823841709473,6.386267964533723);
\draw [color=xfqqff] (-2.2594823841709473,6.386267964533723)-- (-2.2415964568948703,5.724488655318839);
\draw [color=xfqqff] (-2.2415964568948703,5.724488655318839)-- (-2.6708587115207436,5.724488655318839);
\draw [color=xfqqff] (-2.6708587115207436,5.724488655318839)-- (-2.6708587115207436,6.958617637368217);
\draw [line width=5.2pt,color=ffcqcb] (-2.2415964568948703,5.724488655318839)-- (-1.8481060568211545,5.152138982484345);
\draw (-1.8481060568211545,5.152138982484345)-- (-1.3115282385388163,5.152138982484345);
\draw [line width=5.2pt,color=ffcqcb] (-1.3115282385388163,5.152138982484345)-- (-0.5603192929435452,4.2220707641282935);
\draw [line width=5.2pt,color=ffcqcb] (-0.5900178651380301,8.161961715570579)-- (0.25702791899225935,7.078827808345808);
\draw [color=xfqqff] (0.25702791899225935,7.078827808345808)-- (0.686290173618133,6.488592208235236);
\draw [color=xfqqff] (0.686290173618133,6.488592208235236)-- (0.25702791899225935,6.488592208235236);
\draw [color=xfqqff] (0.25702791899225935,6.488592208235236)-- (0.25702791899225935,7.078827808345808);
\draw [line width=5.2pt,color=ffcqcb] (0.686290173618133,6.488592208235236)-- (1.0618946464157715,5.9877862445050525);
\draw [color=xfqqff] (1.0618946464157715,5.9877862445050525)-- (1.4732709737655632,5.415436571670558);
\draw [color=xfqqff] (1.4732709737655632,5.415436571670558)-- (1.491156901041641,4.753657262455674);
\draw [color=xfqqff] (1.491156901041641,4.753657262455674)-- (1.0618946464157715,4.753657262455674);
\draw [color=xfqqff] (1.0618946464157715,4.753657262455674)-- (1.0618946464157715,5.9877862445050525);
\draw [line width=5.2pt,color=ffcqcb] (1.491156901041641,4.753657262455674)-- (1.8846473011153573,4.1813075896211815);
\draw (1.8846473011153573,4.1813075896211815)-- (2.42122511939769,4.1813075896211815);
\draw [color=xfqqff] (2.599087002017442,8.254552666897432)-- (3.028349256643314,7.664317066786861);
\draw [color=xfqqff] (3.028349256643314,7.664317066786861)-- (2.599087002017442,7.664317066786861);
\draw [color=xfqqff] (2.599087002017442,7.664317066786861)-- (2.599087002017442,8.254552666897432);
\draw [line width=5.2pt,color=ffcqcb] (3.028349256643314,7.664317066786861)-- (4.012075256827602,6.287100666528859);
\draw [color=xfqqff] (4.012075256827602,6.287100666528859)-- (4.4234515841774,5.714750993694366);
\draw [color=xfqqff] (4.4234515841774,5.714750993694366)-- (4.441337511453479,5.052971684479483);
\draw [color=xfqqff] (4.441337511453479,5.052971684479483)-- (4.012075256827602,5.052971684479483);
\draw [color=xfqqff] (4.012075256827602,5.052971684479483)-- (4.012075256827602,6.287100666528859);
\draw [line width=5.2pt,color=ffcqcb] (4.441337511453479,5.052971684479483)-- (4.9958012570118875,4.248104957055975);
\draw (4.9958012570118875,4.248104957055975)-- (5.532379075294227,4.248104957055975);
\begin{scriptsize}
\draw [fill=qqqqff] (-6.856165694122984,7.942343637552505) circle (1.5pt);
\draw [fill=qqqqff] (-6.462675294049269,7.405765819270167) circle (1.5pt);
\draw [fill=qqqqff] (-6.033413039423397,6.815530219159594) circle (1.5pt);
\draw [fill=qqqqff] (-6.462675294049269,6.815530219159594) circle (1.5pt);
\draw [fill=qqqqff] (-5.657808566625761,6.314724255429412) circle (1.5pt);
\draw [fill=qqqqff] (-5.246432239275968,5.742374582594918) circle (1.5pt);
\draw [fill=qqqqff] (-5.2285463119998905,5.080595273380034) circle (1.5pt);
\draw [fill=qqqqff] (-5.657808566625761,5.080595273380034) circle (1.5pt);
\draw [fill=qqqqff] (-4.835055911926176,4.50824560054554) circle (1.5pt);
\draw [fill=qqqqff] (-4.298478093643838,4.50824560054554) circle (1.5pt);
\draw [fill=qqqqff] (-3.8871017662940446,3.989553709539279) circle (1.5pt);
\draw [fill=qqqqff] (-3.4757254389442522,8.049659201208971) circle (1.5pt);
\draw [fill=qqqqff] (-3.046463184318381,7.4594236010984) circle (1.5pt);
\draw [fill=qqqqff] (-3.4757254389442522,7.4594236010984) circle (1.5pt);
\draw [fill=qqqqff] (-2.6708587115207436,6.958617637368217) circle (1.5pt);
\draw [fill=qqqqff] (-2.2594823841709473,6.386267964533723) circle (1.5pt);
\draw [fill=qqqqff] (-2.2415964568948703,5.724488655318839) circle (1.5pt);
\draw [fill=qqqqff] (-2.6708587115207436,5.724488655318839) circle (1.5pt);
\draw [fill=qqqqff] (-1.8481060568211545,5.152138982484345) circle (1.5pt);
\draw [fill=qqqqff] (-1.3115282385388163,5.152138982484345) circle (1.5pt);
\draw [fill=qqqqff] (-0.5603192929435452,4.2220707641282935) circle (1.5pt);
\draw [fill=qqqqff] (-0.5900178651380301,8.161961715570579) circle (1.5pt);
\draw [fill=qqqqff] (0.25702791899225935,7.078827808345808) circle (1.5pt);
\draw [fill=qqqqff] (0.686290173618133,6.488592208235236) circle (1.5pt);
\draw [fill=qqqqff] (0.25702791899225935,6.488592208235236) circle (1.5pt);
\draw [fill=qqqqff] (1.0618946464157715,5.9877862445050525) circle (1.5pt);
\draw [fill=qqqqff] (1.4732709737655632,5.415436571670558) circle (1.5pt);
\draw [fill=qqqqff] (1.491156901041641,4.753657262455674) circle (1.5pt);
\draw [fill=qqqqff] (1.0618946464157715,4.753657262455674) circle (1.5pt);
\draw [fill=qqqqff] (1.8846473011153573,4.1813075896211815) circle (1.5pt);
\draw [fill=qqqqff] (2.42122511939769,4.1813075896211815) circle (1.5pt);
\draw [fill=qqqqff] (2.599087002017442,8.254552666897432) circle (1.5pt);
\draw [fill=qqqqff] (3.028349256643314,7.664317066786861) circle (1.5pt);
\draw [fill=qqqqff] (2.599087002017442,7.664317066786861) circle (1.5pt);
\draw [fill=qqqqff] (4.012075256827602,6.287100666528859) circle (1.5pt);
\draw [fill=qqqqff] (4.4234515841774,5.714750993694366) circle (1.5pt);
\draw [fill=qqqqff] (4.441337511453479,5.052971684479483) circle (1.5pt);
\draw [fill=qqqqff] (4.012075256827602,5.052971684479483) circle (1.5pt);
\draw [fill=qqqqff] (4.9958012570118875,4.248104957055975) circle (1.5pt);
\draw [fill=qqqqff] (5.532379075294227,4.248104957055975) circle (1.5pt);
\end{scriptsize}
\end{tikzpicture}

Donc pour obtenir $i$ pentagone dégénéré le nombre de façons de partager les $n-d/r$ pas en $i$ parties est égale à:

\[ \binom{n-d/r-1}{i-1}+2\binom{n-d/r-1}{i}+\binom{n-d/r-1}{i+1}=\binom{n-d/r+1}{i+1} \]
\\

En utilisant deux fois la règle de Pascal.
\\

\begin{note}  $\mu \vDash (n,k)$ est la notation pour une composition, $\mu=\mu_1, \cdots, \mu_k$, de $n$ avec des zéros en $k$ parties, où $\mu_i \in \mathbb{N}$.
\end{note}
\vspace{-35pt}
\begin{ex*} $[0,3],[3,0],[1,2] \vDash (3,2)$, car $0+3=3,~3+0=3,~1+2=3$ sont des sommes égales à trois contenant toutes exactement deux termes.
\end{ex*}
Comme la séparation par des diagonales assure que chaque partie prend ses valeurs pour l'aire gauche de ligne dans des ensembles disjoints. Le nombre de permutations d'aire gauche pour $i$ pentagones ayant au total $d$ pas vers le bas et $d/r$ pas vers la droite est données par:
\[\Pag{_i}= \sum_{\alpha  \vDash (d,i) \atop \beta \vDash (d/r,i)} \binom{d}{\alpha}\prod_{k=1}^{i} \PP(v) \chi(_{\alpha_k\not=0 ~ou~ \beta_k \not=0})\chi(_{\sum_{j=1}^{k} \alpha_j-r\beta_i \geq 0}), \]
où $v=(\alpha_k,\beta_k, \min(\beta_k,\sum_{j=1}^{k-1} \frac{\alpha_j}{r}-\beta_i),r,\sum_{j=1}^{k} \alpha_j-r\beta_i)$ et $\chi$ est le booléen donnant $1$ si c'est vrai et $0$ si c'est faux. 
\\

Notons que dans $\PP$ la valeur de $p$ est le minimum de $\beta_k,\sum_{j=1}^{k-1} \frac{\alpha_j}{r}-\beta_i$, car la hauteur de la partie qui précède est $\sum_{j=1}^{k-1} \alpha_j$ et comme les chemins dans le pentagone ne doivent pas dépasser la diagonale principale qui a une pente de $1/r$ nous devons avoir au plus $\sum_{j=1}^{k-1} \frac{\alpha_j}{r}-\beta_i$, mais comme il est impossible d'avoir des pas vers la gauche $p$ doit être d'au plus $\beta_k$. Si la somme des $\beta_i$ jusqu'à $k$ fois $r$ est strictement plus petite que la somme des $\alpha_i$ jusqu'à $k$, nous avons que $\beta_k$ ne croise pas la diagonale principale donc la valeur de $q$ qui est la hauteur à laquelle on a plus de diagonales principales est bien $\sum_{j=1}^{k} \alpha_j-r\beta_i$. De plus $\chi(_{\sum_{j=1}^{k} \alpha_j-r\beta_i \geq 0})$ nous assure que la valeur de $q$ soit positive, car le contraire serait absurde et $\chi(_{\alpha_k\not=0 ~ou~ \beta_k \not=0})$ nous assure que nous avons bien $i$ pentagones. 

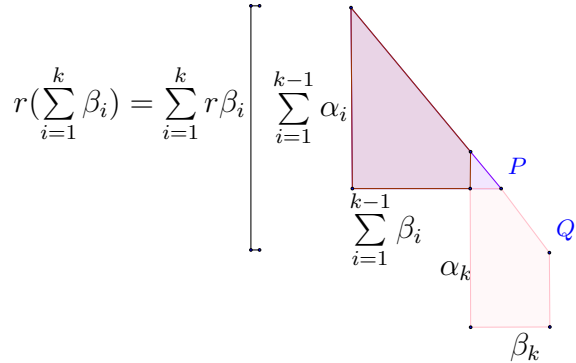
\definecolor{zzttqq}{rgb}{0.6,0.2,0.0}
\definecolor{ffcqcb}{rgb}{1.0,0.7529411764705882,0.796078431372549}
\definecolor{xfqqff}{rgb}{0.4980392156862745,0.0,1.0}
\definecolor{qqqqff}{rgb}{0.0,0.0,1.0}
\begin{tikzpicture}[line cap=round,line join=round,>=triangle 45,x=1.5cm,y=1.5cm]
\clip(-10.864554515294726,4.326178076425408) rectangle (-4.495802191763074,8.245571407672825);
\fill[color=xfqqff,fill=xfqqff,fill opacity=0.1] (-5.660380161393642,6.308052379211532) -- (-6.97894149814047,6.308052379211525) -- (-6.991038391138147,7.912906850236707) -- cycle;
\fill[color=ffcqcb,fill=ffcqcb,fill opacity=0.1] (-5.660380161393642,6.308052379211532) -- (-5.934576402674327,6.308052379211531) -- (-5.930544105008434,5.07820159111434) -- (-5.2285463119998905,5.080595273380034) -- (-5.232956608809045,5.739498408320701) -- cycle;
\fill[color=zzttqq,fill=zzttqq,fill opacity=0.1] (-5.934576402674327,6.308052379211531) -- (-5.931266860544336,6.634758155762972) -- (-6.991038391138147,7.912906850236707) -- (-6.97894149814047,6.308052379211525) -- cycle;
\draw [color=xfqqff] (-5.660380161393642,6.308052379211532)-- (-6.97894149814047,6.308052379211525);
\draw [color=xfqqff] (-6.97894149814047,6.308052379211525)-- (-6.991038391138147,7.912906850236707);
\draw [color=xfqqff] (-6.991038391138147,7.912906850236707)-- (-5.660380161393642,6.308052379211532);
\draw (-7.776174204065816,7.469858144196774) node[anchor=north west] {$\tiny{\sum\limits_{i=1}^{k-1} \alpha_i}$};
\draw (-7.098393405921858,6.38959630445134) node[anchor=north west] {$\tiny{\sum\limits_{i=1}^{k-1} \beta_i}$};
\draw (-6.301363788050991,5.759258631111832) node[anchor=north west] {$\tiny{\alpha_k}$};
\draw (-5.678131566176417,5.12230484683504) node[anchor=north west] {$\tiny{\beta_k}$};
\draw (-7.801530221982523,7.929036040900277)-- (-7.8781438776344785,7.929036040900277);
\draw (-7.8781438776344785,7.929036040900277)-- (-7.878143877634479,5.763692194316049);
\draw (-7.878143877634479,5.763692194316049)-- (-7.801530221982523,5.763692194316049);
\draw (-10.0787654182246,7.510685158063935) node[anchor=north west] {$\tiny{r(\sum\limits_{i=1}^{k} \beta_i)=\sum\limits_{i=1}^{k} r \beta_i}$};
\draw [color=ffcqcb] (-5.660380161393642,6.308052379211532)-- (-5.934576402674327,6.308052379211531);
\draw [color=ffcqcb] (-5.934576402674327,6.308052379211531)-- (-5.930544105008434,5.07820159111434);
\draw [color=ffcqcb] (-5.930544105008434,5.07820159111434)-- (-5.2285463119998905,5.080595273380034);
\draw [color=ffcqcb] (-5.2285463119998905,5.080595273380034)-- (-5.232956608809045,5.739498408320701);
\draw [color=ffcqcb] (-5.232956608809045,5.739498408320701)-- (-5.660380161393642,6.308052379211532);
\draw [color=zzttqq] (-5.934576402674327,6.308052379211531)-- (-5.931266860544336,6.634758155762972);
\draw [color=zzttqq] (-5.931266860544336,6.634758155762972)-- (-6.991038391138147,7.912906850236707);
\draw [color=zzttqq] (-6.991038391138147,7.912906850236707)-- (-6.97894149814047,6.308052379211525);
\draw [color=zzttqq] (-6.97894149814047,6.308052379211525)-- (-5.934576402674327,6.308052379211531);
\begin{scriptsize}
\draw [fill=qqqqff] (-5.660380161393642,6.308052379211532) circle (0.5pt);
\draw[color=qqqqff] (-5.516477538442095,6.510423318318501) node {$P$};
\draw [fill=qqqqff] (-5.232956608809045,5.739498408320701) circle (0.5pt);
\draw[color=qqqqff] (-5.087793892836907,5.938845124178252) node {$Q$};
\draw [fill=qqqqff] (-5.2285463119998905,5.080595273380034) circle (0.5pt);
\draw [fill=qqqqff] (-6.97894149814047,6.308052379211525) circle (0.5pt);
\draw [fill=qqqqff] (-6.991038391138147,7.912906850236707) circle (0.5pt);
\draw [fill=qqqqff] (-7.801530221982523,7.929036040900277) circle (0.5pt);
\draw [fill=qqqqff] (-7.8781438776344785,7.929036040900277) circle (0.5pt);
\draw [fill=qqqqff] (-7.878143877634479,5.763692194316049) circle (0.5pt);
\draw [fill=qqqqff] (-7.801530221982523,5.763692194316049) circle (0.5pt);
\draw [fill=qqqqff] (-5.934576402674327,6.308052379211531) circle (0.5pt);
\draw [fill=qqqqff] (-5.930544105008434,5.07820159111434) circle (0.5pt);
\draw [fill=qqqqff] (-5.931266860544336,6.634758155762972) circle (0.5pt);
\end{scriptsize}
\end{tikzpicture}

\begin{prop} Le nombre de fonctions de stationnement associé aux $r$-FSchröder est donné par la formule:

\begin{align*}  \Ps{_{n,d}^r}&=\binom{nr}{d,\underbrace{r,\cdots,r}_\textrm{n-d/r fois}}\sum_{i=0}^{d+d/r} \binom{n-d/r+1}{i+1}\Pag{_i},
\\ = \binom{nr}{d,\underbrace{r,\cdots,r}_\textrm{n-d/r fois}}&\sum_{i=0}^{d+d/r} \binom{n-d/r+1}{i+1}\sum_{\alpha  \vDash (d,i) \atop \beta \vDash (d/r,i)} \binom{d}{\alpha}\prod_{k=1}^{i} \PP(v) \chi(_{\alpha_k\not=0 ~ou~ \beta_k \not=0})\chi(_{\sum_{j=1}^{k} \alpha_j-r\beta_i \geq 0}), \end{align*}
\end{prop}

où $v=(\alpha_k,\beta_k, \min(\beta_k,\sum\limits_{j=1}^{k-1} \frac{\alpha_j}{r}-\beta_i),r,\sum_{j=1}^{k} \alpha_j-r\beta_i)$ 

\begin{pre} Par construction $\F$
\end{pre}

\section{Fonctions de stationnement sur les chemins 
\\ sans diagonale incluse dans un hexagone.}

 \begin{de} Un \begin{bf} chemin de $r$-Schröder avec contrainte (respectivement les $r$-FSchröder avec contrainte )\end{bf} est un chemin de $r$-Schröder (respectivement un chemin de $r$-FSchröder), n'ayant aucun pas sous la droite de pente $-r$ passant par $(0,n-h)$ (respectivement sous la droite de pente $-1/r$ passant par $(0,nr-h)$)\end{de}

Dans un chemin de Schröder avec contrainte la forme entre deux diagonales non consécutive est un hexagone possiblement dégénéré. Afin d'établir une formule comptant les chemins de stationnement pour les $r$-FSchröder ayant une contrainte $h$, nous allons dans cette section établir une formule comptant le nombre de chemins de stationnements sans diagonale sont dans un hexagone. 
\\

Soit un hexagone défini par les points:

\[ (0,a), (p,a), (b,q), (b,0),(c,0), (0,s),\]
\[  a, b, s, q, r \in \mathbb{N}, p, c \in \mathbb{R}^+, \frac{1}{r}=\frac{a-q}{b-p}=\frac{s}{c}, a \text{ ou } b \text{ non nul et si } a=c,\text{ alors } a=0 \text{ ou } a\not= p\]

\newpage

Par abus nous nommerons cette forme hexagone, mais plusieurs cas sont possibles:
\\Cas 1, la forme est un pentagone tel que défini à la section (\ref{pentagone}{}).
\\Cas 2, $(b-c)+p \leq b$.
\\Cas 3, $(b-c)+p > b$.
\\

Commençons par le cas 2. Si $(b-c)+p \leq b$, les chemins dans l'hexagone doivent tous avoir un pas vers la droite débutant par $(p,i), i \in \{ 0, \cdots, a-1\}$. Cette colonne divise l'hexagone en deux pentagones, l'un est tel que défini à la section (\ref{pentagone}{}) et l'autre est une rotation d'un angle $\pi$ d'un pentagone défini à la section (\ref{pentagone}{}). Comme les rotations ne changent rien au nombre de fonctions de stationnement contenues dans une forme et que l'ensemble des valeurs possible pour l'aire gauche d'une ligne avant un pas vers la droite est disjoint de l'ensemble des valeurs possible pour l'aire gauche d'une ligne après un pas vers la droite, nous pouvons faire le produit de combien de fonctions de stationnement sont possible avant le pas droit par le nombre de possibilités après le pas et par comment permuter les deux ensembles. 
\\ 

De plus, le cas 2 inclu le cas 1 puisque la partie avant la colonne $p$ est alors un segment ou un rectangle qui ont été définie comme un pentagone à la section (\ref{pentagone}).
\\

De façon similaire au cas 2, nous pouvons dans le cas 3 séparer l'hexagone en un pentagone et un hexagone plus petit par un pas vers la droite dans la colonne $p$. Puisque les chemins dans l'hexagone doivent tous avoir un pas vers la droite débutant par $(p,i)$ et que celui-ci ne peut pas être sous la diagonale de la contrainte nous devons avoir  $i \in \{ \max(0,s-pr), \cdots, a-1\}$. Comme ici les pentagones peuvent être vus comme des hexagones, le cas 2 est inclus dans le cas 3.
\\

Nous avons alors la formule (À vérifier) :

\begin{equation*} \hexa(a, b, p, r, q, s)= \sum_{i=\max(0,s-pr)}^{q+(b-\lfloor p \rfloor-1)r} \PP (v)\hexa (w)\binom{a}{i}
\end{equation*} 

Où $v=(a-i,\lfloor p \rfloor,\lfloor p \rfloor-\max\left(\frac{s-i}{r},0\right),r,\min(a-s,a-i))$
\\

et $w=(i,b-\lfloor p \rfloor-1,p- \lfloor p \rfloor +\frac{a-i}{r},r,\min(i,q),i)$.

\newpage

\section{Fonctions de stationnement sur les chemins 
\\ $r$-FSchröder ayant une contrainte}
  
Tout comme à la section (\ref{id}{}) il suffit maintenant de multiplier les parties non diagonales au  multinôme: 

\[ \binom{nr}{d,\underbrace{r,\cdots,r}_\textrm{n-d/r fois}} \]
\\ 

Et de la même façon qu'à la section (\ref{id}{}) $i$ hexagone dégénérée le nombre de façons de partager les $n-d/r$ pas en $i$ parties est égal à:

\[ \binom{n-d/r+1}{i+1} \] 
\\

De la même façon qu'à la section (\ref{id}{}) nous trouvons le nombre de permutations d'aire gauche pour $i$ hexagones ayant au total $d$ pas vers le bas et $d/r$ pas vers la droite et nous trouvons les valeurs avec lesquelles évaluer $\hexa$ pour $a, b, p, r, q$. Il suffit alors de remarquer que la valeur pour $s$ est le minimum entre $h$ et la somme de $1$ à $k$ des $\alpha_j$ D'où la formule (à vérifier):

\begin{align*}  \Cs{_{n,d}^r}&=\binom{nr}{d,\underbrace{r,\cdots,r}_\textrm{n-d/r fois}}\sum_{i=0}^{d+d/r} \binom{n-d/r+1}{i+1}\sum_{\alpha  \vDash (d,i) \atop \beta \vDash (d/r,i)} \binom{d}{\alpha}\prod_{k=1}^{i} \hexa(v) \chi(_{\alpha_k\not=0 ~ou~ \beta_k \not=0})\chi(_{\sum_{j=1}^{k} \alpha_j-r\beta_i \geq 0}), 
\end{align*}

où $v=(\alpha_k,\beta_k,\min(\beta_k,\sum\limits_{j=1}^{k-1} \frac{\alpha_j}{r}-\beta_i),r,\sum_{j=1}^{k} \alpha_j-r\beta_i,\min(h,\sum_{j=1}^{k} \alpha_j)).$
 
\newpage

\section*{\label{Annexe 1}{Annexe 1}}

\vspace{20pt} 
\small{Soit $\mathcal{P} \in \p(0012)$. Rappelons que $\A(\alpha_1)=0,~\A(\alpha_2)=0,~\A(\alpha_3)=1,~\A(\alpha_4)=2.$}

\vspace{50pt} 

$\begin{matrix}
[\mathcal{P}(1),\mathcal{P}(2),\mathcal{P}(3),\mathcal{P}(4)] ,&\hspace{50pt} q^w(\sum_{i=1}^4 q^{\A(\alpha_i)-\mathcal{P}(i)}), \hspace{50pt} &w,
\\ &&
\\ &&
\\  [0,0,1,2],& 4,           & 0,
\\ &&
\\  [0,0,2,1],& 1+2q+q^2,    &1,
\\ &&
\\  [0,2,1,0],& 1+2q^2+q^4,  &2,
\\ &&
\\ [0,2,0,1],& 1+q^2+2q^3,  &2,
\\ &&
\\ [0,1,2,0], & 2+q+q^3,     &1,
\\ &&
\\ [0,1,0,2], & 1+2q+q^2,    &1,
\\ &&
\\ [2,0,1,0], & 1+2q^2+q^4,  &2,
\\ &&
\\ [2,1,0,0], & 1+q+q^3+q^4, &2,
\\ &&
\\ [1,0,2,0], & 2+q+q^3,    &1,
\\ &&
\\ [2,0,0,1], & 1+q^2+2q^3,    &2,
\\ &&
\\ [1,2,0,0], & 1+q+q^3+q^4,  &2,
\\ &&
\\ [1,0,0,2], & 1+2q+q^2,    &1.
\end{matrix}$

\end{document}